\documentclass[final,5p,times]{elsarticle}

\usepackage[utf8]{inputenc}
\usepackage[T1]{fontenc}
\usepackage{amsmath,amssymb,amsfonts}
\usepackage{microtype}
\usepackage{graphicx}
\usepackage{url}
\usepackage[hidelinks]{hyperref}
\usepackage{booktabs}
\usepackage{array}
\usepackage{xcolor}
\usepackage{siunitx}
\usepackage{textcomp}
\DeclareSIPrefix\micro{\text{\textmu}}{-6}
\usepackage{multirow}
\usepackage{acro}
\usepackage{placeins}

\DeclareAcronym{cmp}{short=CMP, long=chemical-mechanical planarisation}
\DeclareAcronym{kvs}{short=KVS, long=K\'arm\'an vortex street}
\DeclareAcronym{wss}{short=WSS, long=wall shear stress}
\DeclareAcronym{cumwss}{short=cum-WSS, long=cumulative wall shear stress}
\DeclareAcronym{pod}{short=POD, long=proper orthogonal decomposition}
\DeclareAcronym{fem}{short=FEM, long=finite element method}
\DeclareAcronym{dns}{short=DNS, long=direct numerical simulation}
\DeclareAcronym{rmse}{short=RMSE, long=root-mean-square error}
\DeclareAcronym{dtw}{short=DTW, long=dynamic time warping}
\DeclareAcronym{pde}{short=PDE, long=partial differential equation}
\DeclareAcronym{rom}{short=ROM, long=reduced-order model}
\DeclareAcronym{dmdc}{short=DMDc, long=dynamic mode decomposition with control}
\DeclareAcronym{mrr}{short=MRR, long=material removal rate}
\DeclareAcronym{wiwnu}{short=WIWNU, long=within-wafer non-uniformity}
\DeclareAcronym{mlp}{short=MLP, long=multilayer perceptron}

\newcommand{\Reyn}{\mathrm{Re}}
\newcommand{\St}{\mathrm{St}}
\newcommand{\uvec}{\mathbf{u}}
\newcommand{\xvec}{\mathbf{x}}
\newcommand{\latent}{\mathbf{z}}
\newcommand{\fftfreq}{E_{\mathrm{St}}}
\newcommand{\rmsew}{\mathrm{RMSE}_{\mathrm{w}}}

\newcommand{\rowKvsInterpFieldDON}{$0.664_{\pm 0.003}$ & $0.717_{\pm 0.002}$ & $0.003_{\pm 0.001}$ & $0.009_{\pm 0.002}$ & $0.004_{\pm 0.002}$ & $0.116_{\pm 0.024}$ & $0.565_{\pm 0.032}$}
\newcommand{\rowKvsInterpSFourD}{$0.340_{\pm 0.034}$ & $0.265_{\pm 0.006}$ & $0.022_{\pm 0.005}$ & $0.021_{\pm 0.002}$ & $0.903_{\pm 0.000}$ & $0.096_{\pm 0.007}$ & $0.052_{\pm 0.004}$}
\newcommand{\rowKvsInterpMamba}{$0.326_{\pm 0.028}$ & $0.262_{\pm 0.006}$ & $0.012_{\pm 0.003}$ & $0.018_{\pm 0.001}$ & $\mathbf{0.907}_{\pm 0.002}$ & $0.093_{\pm 0.013}$ & $0.058_{\pm 0.007}$}
\newcommand{\rowKvsInterpGRU}{$0.440_{\pm 0.043}$ & $0.291_{\pm 0.009}$ & $0.003_{\pm 0.001}$ & $0.020_{\pm 0.001}$ & $0.898_{\pm 0.007}$ & $0.086_{\pm 0.004}$ & $0.042_{\pm 0.003}$}
\newcommand{\rowKvsInterpARMLP}{$0.942_{\pm 0.020}$ & $0.438_{\pm 0.018}$ & $0.013_{\pm 0.008}$ & $0.018_{\pm 0.002}$ & $0.839_{\pm 0.094}$ & $0.388_{\pm 0.029}$ & $0.165_{\pm 0.023}$}
\newcommand{\rowKvsInterpDMDc}{$1.379$ & $0.766$ & $0.042$ & $0.024$ & $0.000$ & --- & $0.891$}
\newcommand{\rowKvsInterpOneshotGRU}{$0.699_{\pm 0.005}$ & $0.718_{\pm 0.012}$ & $0.011_{\pm 0.001}$ & $0.012_{\pm 0.001}$ & $0.032_{\pm 0.021}$ & $0.154_{\pm 0.015}$ & $0.601_{\pm 0.056}$}
\newcommand{\rowKvsInterpOneshotSFourD}{$0.717_{\pm 0.001}$ & $0.740_{\pm 0.000}$ & $0.016_{\pm 0.002}$ & $0.017_{\pm 0.002}$ & $0.000_{\pm 0.000}$ & --- & $0.889_{\pm 0.002}$}
\newcommand{\rowKvsExtrapFieldDON}{$1.097_{\pm 0.005}$ & $1.379_{\pm 0.001}$ & $0.008_{\pm 0.003}$ & $0.002_{\pm 0.000}$ & $0.002_{\pm 0.000}$ & $0.066_{\pm 0.004}$ & $0.311_{\pm 0.105}$}
\newcommand{\rowKvsExtrapSFourD}{$1.292_{\pm 0.024}$ & $0.776_{\pm 0.024}$ & $0.017_{\pm 0.003}$ & $0.037_{\pm 0.011}$ & $\mathbf{0.748}_{\pm 0.014}$ & $0.251_{\pm 0.022}$ & $0.015_{\pm 0.005}$}
\newcommand{\rowKvsExtrapMamba}{$1.328_{\pm 0.032}$ & $0.832_{\pm 0.027}$ & $0.014_{\pm 0.006}$ & $0.031_{\pm 0.017}$ & $0.604_{\pm 0.073}$ & $0.248_{\pm 0.045}$ & $0.017_{\pm 0.002}$}
\newcommand{\rowKvsExtrapGRU}{$1.357_{\pm 0.046}$ & $0.835_{\pm 0.049}$ & $0.009_{\pm 0.001}$ & $0.026_{\pm 0.007}$ & $0.616_{\pm 0.103}$ & $0.258_{\pm 0.012}$ & $0.019_{\pm 0.006}$}
\newcommand{\rowKvsExtrapARMLP}{$1.388_{\pm 0.031}$ & $0.844_{\pm 0.047}$ & $0.016_{\pm 0.002}$ & $0.015_{\pm 0.004}$ & $0.676_{\pm 0.031}$ & $0.355_{\pm 0.171}$ & $0.042_{\pm 0.032}$}
\newcommand{\rowKvsExtrapDMDc}{$2.122$ & $1.454$ & $0.003$ & $0.004$ & $0.000$ & --- & $0.959$}
\newcommand{\rowKvsExtrapOneshotGRU}{$1.224_{\pm 0.011}$ & $1.414_{\pm 0.016}$ & $0.027_{\pm 0.011}$ & $0.011_{\pm 0.004}$ & $0.171_{\pm 0.006}$ & $0.301_{\pm 0.099}$ & $0.717_{\pm 0.014}$}
\newcommand{\rowKvsExtrapOneshotSFourD}{$1.142_{\pm 0.001}$ & $1.414_{\pm 0.001}$ & $0.024_{\pm 0.003}$ & $0.005_{\pm 0.001}$ & $0.000_{\pm 0.000}$ & --- & $0.965_{\pm 0.003}$}
\newcommand{\rowCmpInterpFieldDON}{$0.235_{\pm 0.004}$ & $\phantom{-}0.999_{\pm 0.000}$ & $0.008_{\pm 0.003}$ & $\mathbf{0.027}_{\pm 0.001}$ & $\phantom{-}0.001_{\pm 0.002}$}
\newcommand{\rowCmpInterpSFourD}{$0.393_{\pm 0.063}$ & $\phantom{-}0.929_{\pm 0.005}$ & $0.010_{\pm 0.003}$ & $0.130_{\pm 0.002}$ & $\phantom{-}0.045_{\pm 0.004}$}
\newcommand{\rowCmpInterpMamba}{$0.538_{\pm 0.246}$ & $\phantom{-}0.903_{\pm 0.012}$ & $0.018_{\pm 0.013}$ & $0.143_{\pm 0.015}$ & $\phantom{-}0.047_{\pm 0.026}$}
\newcommand{\rowCmpInterpGRU}{$0.332_{\pm 0.017}$ & $\phantom{-}0.914_{\pm 0.008}$ & $0.004_{\pm 0.001}$ & $0.167_{\pm 0.021}$ & $-0.009_{\pm 0.008}$}
\newcommand{\rowCmpInterpARMLP}{$0.998_{\pm 0.118}$ & $\phantom{-}0.531_{\pm 0.151}$ & $0.078_{\pm 0.003}$ & $0.419_{\pm 0.023}$ & $-0.027_{\pm 0.019}$}
\newcommand{\rowCmpInterpDMDc}{$2.379$ & $-0.110$ & $0.575$ & $0.796$ & $-0.121$}
\newcommand{\rowCmpInterpOneshotGRU}{$0.326_{\pm 0.014}$ & $\phantom{-}0.968_{\pm 0.011}$ & $0.002_{\pm 0.000}$ & $0.099_{\pm 0.012}$ & $\phantom{-}0.013_{\pm 0.010}$}
\newcommand{\rowCmpInterpOneshotSFourD}{$0.506_{\pm 0.009}$ & $\phantom{-}0.530_{\pm 0.020}$ & $0.029_{\pm 0.004}$ & $0.376_{\pm 0.004}$ & $-0.103_{\pm 0.007}$}
\newcommand{\rowCmpExtrapFieldDON}{$0.306_{\pm 0.014}$ & $\phantom{-}0.999_{\pm 0.000}$ & $0.010_{\pm 0.003}$ & $\mathbf{0.032}_{\pm 0.006}$ & $\phantom{-}0.000_{\pm 0.002}$}
\newcommand{\rowCmpExtrapSFourD}{$0.416_{\pm 0.018}$ & $\phantom{-}0.952_{\pm 0.005}$ & $0.005_{\pm 0.001}$ & $0.123_{\pm 0.005}$ & $\phantom{-}0.010_{\pm 0.005}$}
\newcommand{\rowCmpExtrapMamba}{$0.555_{\pm 0.119}$ & $\phantom{-}0.956_{\pm 0.005}$ & $0.004_{\pm 0.001}$ & $0.120_{\pm 0.004}$ & $\phantom{-}0.016_{\pm 0.009}$}
\newcommand{\rowCmpExtrapGRU}{$0.402_{\pm 0.040}$ & $\phantom{-}0.922_{\pm 0.025}$ & $0.005_{\pm 0.001}$ & $0.151_{\pm 0.020}$ & $-0.004_{\pm 0.007}$}
\newcommand{\rowCmpExtrapARMLP}{$37.7_{\pm 52.1}$ & $\phantom{-}0.763_{\pm 0.070}$ & $0.100_{\pm 0.035}$ & $1.775_{\pm 2.113}$ & $-0.012_{\pm 0.015}$}
\newcommand{\rowCmpExtrapDMDc}{$3.414$ & $\phantom{-}0.303$ & $0.513$ & $0.848$ & $-0.077$}
\newcommand{\rowCmpExtrapOneshotGRU}{$0.372_{\pm 0.003}$ & $\phantom{-}0.974_{\pm 0.002}$ & $0.001_{\pm 0.000}$ & $0.091_{\pm 0.003}$ & $\phantom{-}0.001_{\pm 0.001}$}
\newcommand{\rowCmpExtrapOneshotSFourD}{$0.690_{\pm 0.018}$ & $\phantom{-}0.748_{\pm 0.009}$ & $0.014_{\pm 0.009}$ & $0.355_{\pm 0.008}$ & $-0.083_{\pm 0.002}$}
\newcommand{\vCaeFitH}{1.05}
\newcommand{\vSFourDFitH}{0.54}
\newcommand{\vMambaFitH}{0.64}

\newcommand{\vSFourDTrainH}{1.59}
\newcommand{\vMambaTrainH}{1.69}
\newcommand{\vKvsQstar}{646}
\newcommand{\vKvsOfflineH}{498}
\newcommand{\vKvsTrainQueriesLo}{2.1}
\newcommand{\vKvsTrainQueriesHi}{2.2}
\newcommand{\vFilmIpFactor}{1.4}
\newcommand{\vKvsFieldDonLastP}{0.01}
\newcommand{\vKvsMirrorCloserCases}{63}
\newcommand{\vKvsMirrorNCases}{200}
\newcommand{\vKvsMirrorGainPct}{14}
\newcommand{\vFilmEncFloorInterpCw}{0.065}

\newcommand{\vFilmEncFloorInterpIp}{0.18}

\newcommand{\vKvsPodModesNinetyNine}{36}
\newcommand{\vKvsPodModeOnePct}{53}
\newcommand{\vKvsEncFloorInterpRmse}{0.09}
\newcommand{\vKvsEncFloorInterpRmsew}{0.08}
\newcommand{\vKvsEncFloorInterpP}{0.91}
\newcommand{\vKvsEncFloorExtrapRmse}{0.29}
\newcommand{\vKvsEncFloorExtrapRmsew}{0.33}
\newcommand{\vKvsEncFloorExtrapP}{0.74}
\newcommand{\vKvsOrbitPowerAr}{0.91}
\newcommand{\vKvsOrbitPowerOs}{0.50}
\newcommand{\vKvsOrbitRadiusAr}{0.96}

\newcommand{\vKvsOrbitRadiusOsLate}{0.36}
\newcommand{\vKvsOrbitBandArPct}{92}
\newcommand{\vKvsOrbitBandOsPct}{1}
\newcommand{\vKvsOrbitAliasedHiPct}{85}
\newcommand{\vKvsOrbitAliasedLoPct}{11}
\newcommand{\vKvsOrbitUmaxSplit}{2.5}
\newcommand{\vKvsOrbitTstar}{3}

\newcommand{\vKvsOrbitNall}{782}
\newcommand{\vKvsOrbitInPlaneTruthLate}{0.91}

\newcommand{\vKvsOrbitInPlaneOsLate}{0.25}
\newcommand{\vKvsOrbitOutPlaneTruthLate}{0.42}
\newcommand{\vKvsOrbitOutPlaneOsLate}{0.34}
\newcommand{\vKvsOrbitWraps}{35}
\newcommand{\vKvsLadderLevels}{46}
\newcommand{\vKvsLadderStep}{0.1}
\newcommand{\vKvsLadderTurnsPerLevel}{0.79}
\newcommand{\vKvsLadderNyquistLevels}{72}
\newcommand{\vKvsLadderNyquistFactor}{1.6}
\newcommand{\vCmpGrowDirectRampEnd}{0.26}
\newcommand{\vCmpGrowDirectEnd}{0.29}

\newcommand{\vCmpGrowGruRampEnd}{0.35}

\newcommand{\vCmpGrowMlpRatio}{3.4}

\newcommand{\vCmpGrowDmdcRatio}{2.4}
\newcommand{\vCmpGrowArFloor}{0.043}

\newcommand{\vKvsOrbitSPowerOs}{0.08}
\newcommand{\vKvsOrbitSRadiusAr}{0.98}
\newcommand{\vKvsOrbitSBandArPct}{96}
\newcommand{\vKvsOrbitSBandOsPct}{0.1}

\newcommand{\vKvsOrbitSInPlaneOsLate}{0.21}
\newcommand{\vKvsOrbitSOutPlaneOsLate}{1.60}
\newcommand{\vKvsFftPairsInterp}{1656}
\newcommand{\vKvsFftDroppedInterp}{14}

\newcommand{\vKvsShedPeriodMin}{0.08}
\newcommand{\vKvsShedPeriodMax}{0.67}
\newcommand{\vKvsUmaxCrit}{0.7}

\biboptions{sort&compress}

\journal{Journal of Computational Science}

\title{No Free Lunch in Flow Surrogates under Time-Varying Boundary Conditions:
A Two-Regime Study}

\author[tuc]{Georg Winkler\corref{cor1}}
\author[tuc]{Martin Stoll}
\affiliation[tuc]{organization={Department of Mathematics, Chemnitz University of Technology}, postcode={09111}, city={Chemnitz}, country={Germany}}
\cortext[cor1]{Corresponding author.\newline E-mail: georg.winkler@mathematik.tu-chemnitz.de}

\begin{document}

\begin{frontmatter}

\begin{abstract}
We test whether an architecture that succeeds on a simple flow regime also succeeds on a
richer one, with each trained separately on each regime. We explore two transient
flows under time-varying boundary conditions: the three-dimensional slurry film in
chemical-mechanical planarisation (CMP), central to semiconductor manufacturing, and the
two-dimensional K\'arm\'an vortex street (KVS). Eight surrogate models on one
shared pipeline differ in whether they learn the full field or a latent representation, and in
whether they predict in one shot or step by step. No
single architecture wins both regimes. On the film, a one-shot full-field model reconstructs the
cumulative wall shear stress to $2.7\%$ relative error. On the wake, a latent
autoregressive DeepONet retains $90\%$ of the shedding power that direct and one-shot models
damp to almost zero. The treatment of time decides the outcome. The self-sustained wake calls
for autoregressive feedback and the boundary-driven film for a direct map. Pointwise RMSE hides the damped oscillation on the wake,
compresses the sixfold lead on the film's process target, and picks the damped model
under wake extrapolation. The evaluation scores five physical
questions. Trained surrogates answer queries $10^{3}$ to $10^{4}$ times faster than the
finite-element solver and pay off from the first query beyond the training set on the film and
from the third on the wake. Neither the winning architecture nor its validation holds
across regimes. The choice of surrogate should follow the dynamical character of the
target flow, and its validation should resolve the failure modes.
\end{abstract}


\end{frontmatter}

\section*{Highlights}
\begin{itemize}
  \item One-shot direct models win the film and autoregressive latent models win the wake.
  \item The treatment of time decides the winner, and the representation sets the margin.
  \item RMSE hides wake damping, compresses film margins, and misranks under extrapolation.
  \item The metric suite scores field, structure, invented motion, amplitude, and timing.
  \item The break-even is set by the training simulations, not by the training compute.
\end{itemize}

\section{Introduction}
\label{sec:intro}
\begin{figure}[t]
\centering
  \includegraphics[width=0.95\linewidth]{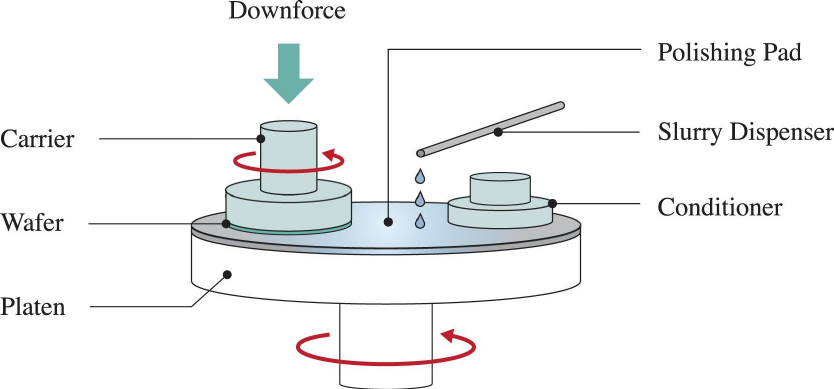}
\caption{Schematic illustration of the \ac{cmp} process: a rotating wafer (held by a carrier under
controlled down-force) is pressed against a slurry-fed, rotating polishing pad. The slurry chemically softens
the surface while the pad and the abrasive particles mechanically abrade it, removing material
to planarise each deposited layer. Schematic reproduced from our
open-access review \cite{winkler2025cmpreview} under the terms of the Creative Commons Attribution
4.0 (CC~BY~4.0) license.}
\label{fig:cmp-schematic}
\end{figure}

\Ac{cmp} is a fundamental process step in the manufacture of
integrated circuits. Modern chips are built as stacks of many patterned layers, and each newly
deposited layer must be polished back to a globally flat surface before the next can be added.
Surface non-uniformities left by imperfect planarisation propagate into defects in the subsequent layers and
reduce manufacturing yield \cite{zantye2004yield}. As its name implies, \ac{cmp} combines a chemical
softening of the wafer surface with mechanical abrasion. A rotating wafer is pressed against a
rotating polishing pad, with an abrasive slurry fed into the interface between them
(Fig.~\ref{fig:cmp-schematic}).

Optimising \ac{cmp} experimentally is hard. The space of process configurations, classically called
the \emph{recipe} in \ac{cmp}, is high-dimensional, spanning many factors such as slurry composition,
operating velocities, and applied pressure. With a wafer consumed per physical characterisation run,
only a handful of configurations can ever be explored.
Compounding this, the coupled chemical-mechanical interactions are strongly non-linear:
first-principles models remain too unreliable to verify or optimise real processes
\cite{krishnan2010cmp}, and conventional statistical process-control methods struggle to capture these
interactions \cite{kim2021r2r}. Process set-up therefore often draws on accumulated
engineering experience to compensate for an incomplete information base, leading to a heightened
dependence on process experts while accepting sub-optimal configurations and the associated yield loss.

Simulating the process and screening configurations in silico rather than in the fab makes
this search tractable, trading scarce silicon for compute and accuracy for speed. The more accurately
the simulation reproduces the process, the more expensive each run becomes, and a screening or
optimisation loop quickly makes simulation time the dominant cost. To speed up this screening, a \emph{surrogate model} trained on simulated data can be
leveraged, predicting the process response and ranking candidates orders of magnitude faster
than the solver, at a further loss of accuracy relative to the simulation it emulates.
Choosing an appropriate architecture is not trivial. The best
choice shifts with the fidelity of the simulation, so a method that excels on a simple process
model may fail on a richer one.

This study compares two contrasting flow regimes. The slurry film in \ac{cmp} is simplified to a
quasi-static 3D Stokes flow. The 2D flow past a cylinder shedding a \ac{kvs} is its counterpart. Its self-sustained
oscillation carries dynamics of a richness that is not to be expected of the simplified
\ac{cmp} flow but could arise in a more realistic model. The pairing targets exactly this contrast in dynamical
character. We refer to the two regimes by their flow features as the \emph{film} and the \emph{wake}. Both are
driven by a gradual ramp in
their boundary conditions. In \ac{cmp}, the ramp-up is necessary to reduce scratch defects at startup, but can introduce
wafer-scale non-uniformities, which are also to be avoided.
The prediction target is therefore the transient flow itself, not a steady snapshot. Predicting
such transients with a surrogate raises three questions. Does a single architecture win both
regimes? Which failure modes emerge, and does a single aggregate error expose them? And when
does the surrogate route pay off against direct simulation? Governing equations, geometries, and the polishing-relevant material-removal target are
deferred to Sec.~\ref{sec:data}.

\paragraph{Contributions}
\begin{enumerate}
  \item \textbf{A break-even cost accounting.} Each query runs three to four orders of
  magnitude faster than the finite-element solver, but the gain only materialises beyond a break-even query
  count set almost entirely by training-set size, $Q^{*}\approx 70$ for the \ac{cmp} film and
  $Q^{*}\approx \vKvsQstar$ for the \ac{kvs} wake. The surrogate is a many-query instrument, not a free lunch
  for individual predictions (Sec.~\ref{sec:results-cost}).
  \item \textbf{A \ac{cmp}-motivated no-free-lunch result.} No single architecture evaluated on a shared pipeline wins both regimes. A one-shot operator that maps the whole trajectory at once wins
  the film, autoregressive feedback that carries the state step by step wins the wake, and neither holds across regimes
  (Sec.~\ref{sec:results}, \ref{sec:discussion}).
  \item \textbf{A problem-specific evaluation suite} that resolves what a pointwise \ac{rmse}
  hides. Limit-cycle damping, pre-onset hallucination, and onset latency each count as mild
  deviations in that error, while the wake's mirror-branch ambiguity under
  extrapolation is over-penalised by it. The suite and the pointwise \ac{rmse} diverge in both
  regimes. On the wake the aggregate hides the damped oscillation, on the film it compresses the
  sixfold lead on the process target into a factor of $\vFilmIpFactor$, and under wake
  extrapolation it inverts the ranking outright
  (Sec.~\ref{sec:metrics},~\ref{sec:discussion}).
  \item \textbf{A validation caution for \ac{cmp} surrogate programmes.} Cheap validation on a
  simplified regime is no evidence of reliability on a richer one: the architecture ranked best on the
  simple \ac{cmp} flow may fail when it is trained and tested on more complex physics
  (Sec.~\ref{sec:discussion}).
\end{enumerate}

The remainder of the paper is structured as follows. Sec.~\ref{sec:related}
situates the study among reduced-order models and operator learning.
Sec.~\ref{sec:data} defines the two regimes and their prediction targets,
Sec.~\ref{sec:methods} describes the surrogate architectures and their training,
and Sec.~\ref{sec:metrics} sets out the evaluation metrics. Sec.~\ref{sec:results}
reports the results for both regimes and the computational cost. Sec.~\ref{sec:discussion} discusses the findings,
Sec.~\ref{sec:limitations} states the limitations together
with the follow-up work they motivate, and Sec.~\ref{sec:conclusion} concludes.

\section{Related Work}
\label{sec:related}

Machine learning has taken on a growing role in \ac{cmp}. The focus lies
on direct applications such as process monitoring \cite{cheng2015time,pan2019scratch}
and the prediction of post-polish quantities
\cite{jebri2016vm,xu2021mrr,deng2021sicmrr,cai2021neurfill}, surveyed in
\cite{winkler2025cmpreview}. Machine learning for the prediction of the slurry flow itself has
so far received little attention.

\ac{cmp} is a closed process, and the flow field between pad and wafer cannot be measured
directly. The physics must therefore come from numerical simulation. Simulating such flows is
well established \cite{ferziger2020cfd}, but the time integration takes many expensive solver
steps. Model reduction emerged early as a remedy that restricts the problem to its dominant
modes. A classical choice is proper orthogonal decomposition with Galerkin projection, which
builds a linear subspace from solution snapshots \cite{berkooz1993pod,benner2015survey}. These
linear reductions are interpretable, but they struggle once the dynamics turn strongly
nonlinear or oscillatory. Data-driven surrogates use a learned nonlinear map instead. Many of
them share one pattern: encode the field into a low-dimensional latent
state \cite{lee2020convae,seidman2023vano}, advance that state with a temporal model, and
decode it back into the original space
\cite{farenga2024latent,oommen2022microstructure,serrano2024aroma}.
This pattern captures dynamics that a linear basis represents only with many modes.
Recent non-intrusive reduced-order models learn the latent evolution recursively. Fu et
al.\ use a self-attention network over past latent codes~\cite{fu2023attention}, Lin et al.\
identify sparse differential equations on the codes~\cite{lin2025discovery}, and Guo et al.\
advance a manifold embedding step by step~\cite{guo2026pmd}. All three demonstrate this on
single trajectories of a cylinder wake at fixed inflow.

The latent surrogates above describe the dynamics of one fixed system and learn a map between
finite-dimensional vectors. Operator learning instead models the solution operator itself, the
map from an input function to the solution function \cite{chen1995universal}. The input can be
a boundary condition or a coefficient field, the output the full flow field. One trained model
then covers a whole class of inputs, at the price that the training data must span that
class \cite{kovachki2024operator}.

Two families realise this idea. Encoder-decoder networks such as PCA-Net
\cite{bhattacharya2020model} and DeepONet compress the input function into a finite number of
coefficients and expand them over learned basis functions. They stay close to the encoding
strategy above. In DeepONet, a branch and a trunk network take these two roles
\cite{lu2021deeponet}, a construction with approximation-theoretic error bounds in infinite
dimensions \cite{lanthaler2022deeponet}. Neural operators, the second family, keep the
function throughout the network. Each layer integrates the input function against a learned
kernel, an operation the Fourier neural operator evaluates as a multiplication in Fourier
space \cite{li2021fno}. Both families have since been refined in many variants
\cite{tripura2023wno,cao2021galerkin,hao2023gnot,rahman2023uno,kontolati2024latent,mandl2025timeintegrated,buitrago2024memno,wang2021improved}.
Across sixteen benchmarks under matched conditions, Lu et al. found no consistent winner
between the two families \cite{lu2022faircomparison}.

Transient problems, as treated in this work, pose the additional question of how to traverse
time. One route predicts the whole trajectory at once. This is a complex learning problem, and
it fixes the prediction to the trained time horizon. The alternative is an autoregressive
rollout, which advances one step at a time and feeds each prediction back as input. The model
becomes more flexible and the learning problem simpler, but small errors accumulate over long
horizons \cite{ye2025recurrent,brandstetter2022message,lippe2023pderefiner}. The instability
has been linked to states that the training data do not
constrain \cite{floryan2024instabilities}. Training the rollout on its own outputs counters
this drift \cite{um2020solver}. The two schemes therefore fail in different ways.

Operator and latent surrogates have been applied to a wide range of evolution problems
\cite{cho2026graphdeeponet,clarkdileoni2021}, and extensions that stabilise long horizons
report good accuracy \cite{wang2022longtime,nayak2025tideeponet}. Across these studies, the varied
input is usually an initial condition, a coefficient field, or a load history.
To the best of our knowledge, a flow driven by a boundary condition that changes during the
transient remains to be addressed. Closest to our setting, Abueidda et al. drive their model
with a time-varying load history \cite{abueidda2025ncde}, though for solid mechanics rather
than a flow.

A general concern is the fair and meaningful comparison of reported performance. McGreivy and
Hakim reviewed 82 studies on machine learning for fluid-related \acp{pde} and found that 79\%
of those claiming to outperform a standard method used weak baselines and that negative
results are systematically underreported, which complicates a sound assessment of the
findings \cite{mcgreivy2024weak}.
Beyond the baselines, the error measure itself can mislead. A low aggregate error can coexist
with unphysical dynamics \cite{mohan2024illusion}, and an error measured against the
solver that generated the training data penalises an emulator that surpasses it
\cite{koehler2025emulator}. A low one-step loss does not imply long-horizon reliability \cite{zou2026uq},
and cost-blind comparisons overstate the practical gain \cite{westermann2026neural}. These findings call for metrics that
resolve where a model fails (Sec.~\ref{sec:metrics}).

\section{Datasets and Problem Setup}
\label{sec:data}
\begin{figure*}[t]
\centering
\includegraphics[width=\textwidth]{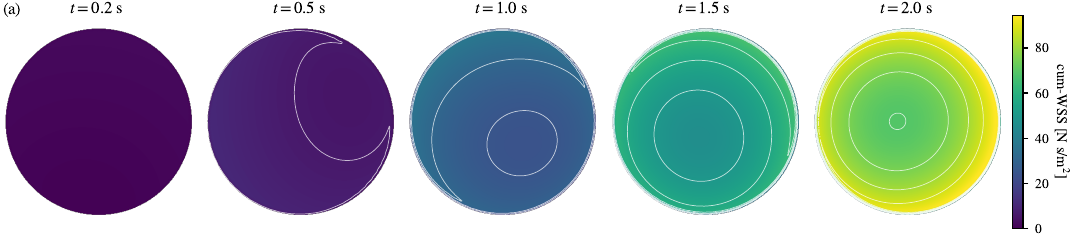}\\[2pt]
\includegraphics[width=\textwidth]{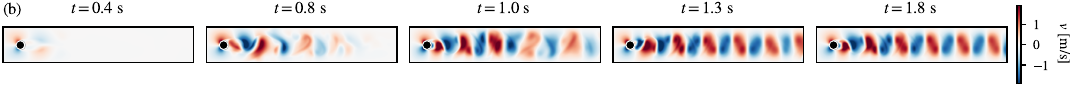}
\caption{Temporal evolution of the ground-truth \ac{fem} solution under the boundary-condition ramp.
(a)~\ac{cmp}: the
\ac{cumwss} on the wafer builds up monotonically over the transient, accumulated in the co-rotating
wafer frame, with white lines marking its iso-levels. This wafer-frame view shows the
\ac{cumwss} received by wafer material points and complements the benchmark quantity of
Eq.~\eqref{eq:cumwss}, which integrates at fixed laboratory positions. (b)~\ac{kvs}: the $v$-velocity field
of a representative case at $\Reyn \approx 167$ develops from near-rest, through the ramp
($t \le 1\,\si{s}$), into self-sustained vortex shedding.}
\label{fig:evolution}
\end{figure*}

\begin{figure*}[t]
\centering
\includegraphics[width=\textwidth]{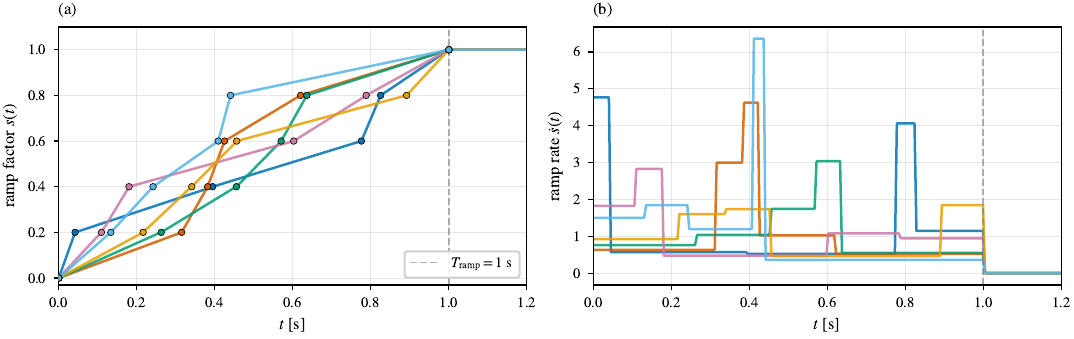}
\caption{The piecewise-linear boundary-condition ramp (representative \ac{kvs} schedules). (a)~Six example
ramps $s(t)$ reach fixed levels $[0.2,0.4,0.6,0.8,1.0]$ at case-dependent knot times $\tau$.
(b)~The corresponding ramp rate $\dot s(t)$ is piecewise-constant between knots. The schedule $\tau$ is
part of the control parameters $\theta$.}
\label{fig:ramps}
\end{figure*}

This study considers the prediction of flow evolutions under time-varying boundary conditions.
The flow of an incompressible fluid can be described by the incompressible Navier--Stokes
equations
\begin{align}
  \partial_t \uvec + (\uvec \cdot \nabla)\uvec
    &= -\tfrac{1}{\rho}\nabla p + \nu \nabla^2 \uvec, &
  \nabla \cdot \uvec &= 0,
  \label{eq:ns}
\end{align}
with velocity $\uvec$, pressure $p$, density $\rho$, and kinematic viscosity $\nu$. In both
datasets, the boundary conditions vary in time following the same construction. The flow is
not switched on at full strength but brought up through a multi-stage, piecewise-linear ramp
$s(t)$ (Fig.~\ref{fig:ramps}), with piecewise-constant rate $\dot s(t)$. The boundary
conditions ramp up over one second, and the case-dependent knot times are drawn as Dirichlet
partitions of the ramp window. The two
datasets instantiate this construction differently. The \ac{kvs} ramps rise over five levels
in steps of $0.2$ at concentration $\alpha = 2$, with $20$ schedules shared across the inflow
grid. The \ac{cmp} ramps rise over
four levels in steps of $0.25$ at $\alpha = 1$, with an individual schedule per case. Within
each regime the schedule family is fixed and shared between training and test. The ramp
schedule is part of the control parameters $\theta$ on
which each surrogate is conditioned. The surrogate must predict the full time-resolved flow fields
this ramp produces (Fig.~\ref{fig:evolution}). Reference solutions come from finite-element
simulations, discretised with Taylor--Hood ($P_2$/$P_1$) elements in DOLFINx
\cite{baratta2023dolfinx}. The time-stepping differs between the two datasets and is noted with each.

\subsection{CMP: Chemical-Mechanical Planarisation}
\label{sec:data-cmp}
Our first dataset is inspired by \ac{cmp}. In simplified terms, the wafer to be polished is
pressed against a slurry-covered, counter-rotating pad, and the slurry flow between the two
influences where material is removed.

To first order, the local material-removal rate is described by Preston's law
\begin{equation}
  \mathrm{MRR} = k_p \, p \, v_{\mathrm{rel}},
  \label{eq:preston}
\end{equation}
with Preston coefficient $k_p$, contact pressure $p$ set by the applied down-force, and
relative pad-wafer velocity $v_{\mathrm{rel}}$, where the influence of the slurry flow enters
\cite{preston1927}. The shear force law generalises this relation and links removal to the
frictional shear force at the wafer \cite{borucki2023shearforce}. Our study considers its
hydrodynamic part, the wall shear stress the slurry exerts on the wafer. The
process target is its time integral at fixed positions in the laboratory frame, the
\ac{cumwss} field
\begin{equation}
  \mathcal{S}(\xvec) = \int_0^T \big\| \tau_w(\xvec, t) \big\|\, dt,
  \qquad
  \tau_w = \mu\, \frac{\uvec_{\parallel,\text{top}} - \uvec_{\parallel,\text{bot}}}{h},
  \label{eq:cumwss}
\end{equation}
where $\tau_w$ is the film-averaged stress across the gap of height $h$ and $\mu = \rho\nu$ the
dynamic viscosity.

\begin{figure*}[t]
\centering
\includegraphics[width=\textwidth]{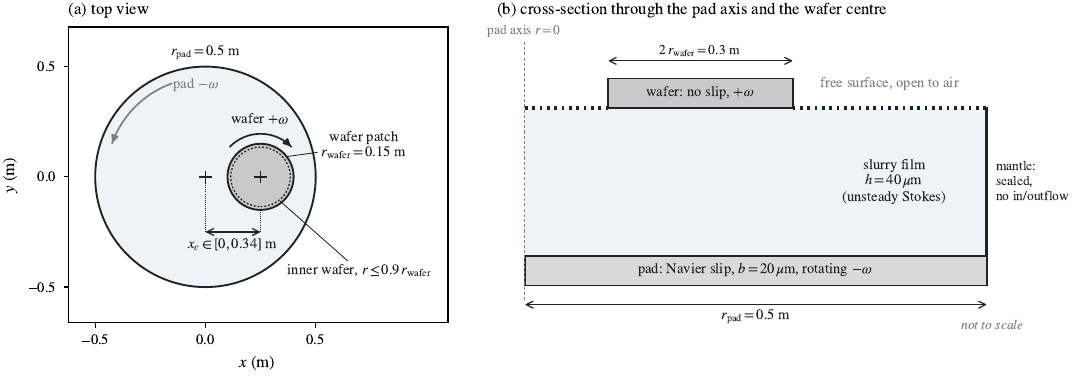}
\caption{\ac{cmp} computational domain. (a)~Domain and wafer placement. The domain is the
whole slurry film over the pad, a flat cylinder of radius $r_{\mathrm{pad}} = 0.5\,\si{\metre}$ and
height $40\,\si{\micro\metre}$. The
wafer enters only as a rotating no-slip patch of radius $r_{\mathrm{wafer}} = 0.15\,\si{\metre}$ on the
upper surface, offset from the pad centre by $x_c$. The dashed circle marks the inner wafer,
$r \le 0.9\,r_{\mathrm{wafer}}$. Sec.~\ref{sec:metrics} reads the film metrics there and the
aggregate error over the whole domain. (b)~Boundary conditions
in the cross-section through the pad axis and the wafer centre, not to scale. The pad below rotates in the opposite direction under a
Navier-slip condition. Outside the wafer, the upper surface is traction-free (free surface to
the air). The
mantle is sealed.}
\label{fig:cmp-domain}
\end{figure*}

We simulate this flow in simplified form. The computational domain is a flat cylinder of
radius $r_{\mathrm{pad}} = 0.5\,\si{\metre}$ and height $40\,\si{\micro\metre}$, the whole
slurry film over the pad and not only the part covered by the wafer
(Fig.~\ref{fig:cmp-domain}). The upper
side carries the wafer, a rotating no-slip disk of radius $r_{\mathrm{wafer}} = 0.15\,\si{\metre}$,
offset from the pad centre by
$x_c \in [0, 0.34]\,\si{\metre}$. The wafer enters as a boundary condition on the upper
surface, not as part of the geometry. A change of its offset therefore only moves that
condition. All cases share one mesh. Outside the wafer the upper surface is traction-free, a
homogeneous Neumann condition that stands in for the free surface to the air without surface
tension or air shear. The
lower side, the pad, carries a Navier-slip condition with slip length
$b = 20\,\si{\micro\metre}$, standing in for a soft pad that yields to the slurry instead of
resolving its texture. The wafer spins at 60 to 120 rpm, the pad at the same rate in the
opposite direction, and both rates follow the ramp. The mantle is sealed by a normal impedance condition, and no slurry enters or leaves the
domain. In the real process, slurry is fed continuously and leaves over the pad edge. The
control parameters $\theta$ are the rotation rate, the offset $x_c$, and the knot times of the
ramp.

In this thin film, the slurry is set in motion almost entirely by shear. The pad below and the
counter-rotating wafer above drag the viscous fluid along at speeds of order $1\,\mathrm{m/s}$.
With the water-like slurry viscosity $\nu = 10^{-6}\,\mathrm{m^2/s}$, the gap Reynolds number is
$Re_h = U h/\nu \approx 40$, and the thin-film aspect ratio $h/r_{\mathrm{wafer}} \approx 3\times10^{-4}$ reduces
the effective inertia to $Re_h\,(h/r_{\mathrm{wafer}}) \approx 10^{-2}$. The convective term is therefore
negligible, and the momentum balance reduces to the unsteady Stokes equations
\begin{equation}
  \rho\,\partial_t \uvec = -\nabla p + \mu \nabla^2 \uvec, \qquad \nabla \cdot \uvec = 0,
  \label{eq:stokes}
\end{equation}
integrated in time by a backward-Euler scheme. We retain the time derivative and resolve
the ramp transient rather than assuming a steady flow.
The viscous term enters in the Laplacian form of Eq.~\eqref{eq:stokes}, whose difference to the
full-stress form is confined to the slip boundary and scales with the gap-to-radius ratio. Viscous relaxation across the gap,
$h^2/\nu \approx 1.6\times10^{-3}\,\si{s}$, is far faster than the one-second ramp, so the
field follows the boundary conditions with negligible lag and the flow turns out quasi-static.

The dataset comprises $100$ cases of $201$ snapshots each, covering $t \in [0,2]\,\si{s}$.
Each snapshot contains the three-dimensional velocity field on $1{,}817{,}169$ finite-element
nodes in $29$ $z$-levels. For training, we subsample this field to $233{,}322$ nodes. The top
and bottom planes keep their full resolution of $62{,}661$ nodes each, since the moving
boundaries drive the flow and the velocity varies most steeply there. The $27$ interior levels are
thinned to $4{,}000$ nodes each. One case takes about $11\,\si{\hour}$ of solver
wall-clock time, detailed in the cost accounting of Sec.~\ref{sec:results-cost}.

For evaluation, \ac{cumwss} is not predicted directly but computed from the predicted velocity
fields. The shear of Eq.~\eqref{eq:cumwss} is evaluated as the difference quotient across the
gap and reported on the $62{,}661$ in-plane positions. The film develops an almost exactly
linear (slip-Couette) profile, so this agrees with the resolved wafer-side wall shear to within
a few percent across the pad land. Localised peaks over under-resolved zones of the FE solution
are excluded by construction, and the same operator is applied to predicted and reference
fields, so the relative film metrics are insensitive to this convention.

\subsection{KVS: K\'arm\'an Vortex Street}
\label{sec:data-kvs}
The second dataset is the flow past a circular cylinder in a two-dimensional channel, which
develops a \ac{kvs}. We choose this flow because it is thoroughly characterised~\cite{williamson1996wake}. The geometry
follows the Sch\"afer--Turek benchmark \cite{schafer1996benchmark}, which spans a startup case
with time-varying inflow (2D-3) and a time-periodic shedding case (2D-2) with a documented
reference frequency at $\Reyn = 100$. After the ramp our flow settles into the same
time-periodic regime.

\begin{figure}[t]
\centering
\includegraphics[width=\linewidth]{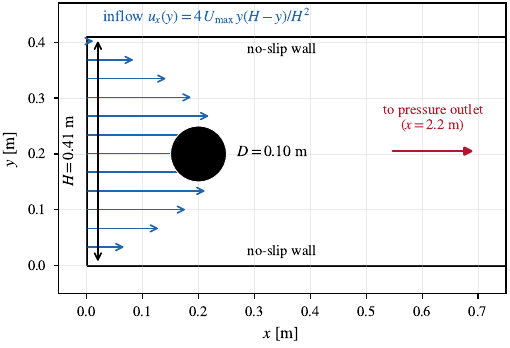}
\caption{\ac{kvs} setup, following the Sch\"afer--Turek benchmark
\cite{schafer1996benchmark}. The parabolic inflow of Eq.~\eqref{eq:inflow} enters from the
left, the channel walls and the cylinder are no-slip, and the channel continues to a pressure
outlet. Only the front region is shown.}
\label{fig:kvs-setup}
\end{figure}

The geometry is taken unchanged from the benchmark and is depicted in
Fig.~\ref{fig:kvs-setup}. A channel of height $H = 0.41$ and length $2.2$ contains a cylinder
of diameter $D = 0.1$, centred at $(0.2, 0.2)$ slightly below the centreline. The flow enters
with a parabolic profile, defined by its peak velocity $U_{\max}$ and scaled in time by the
ramp,
\begin{equation}
  u_x(y, t) = 4\, s(t)\, U_{\max}\, \frac{y\,(H-y)}{H^2}.
  \label{eq:inflow}
\end{equation}
Instead of the benchmark's prescribed inflow history, we vary $U_{\max}$ from $0.5$ to $5.0$
in steps of $0.1$. Every level is crossed with $20$ ramp schedules between $0$ and
$1\,\si{s}$, resulting in $920$ available simulations for modelling. Each case comprises $201$
snapshots of $u$ and $v$ on the $237{,}351$ second-order velocity nodes and of $p$ on the $59{,}691$ first-order pressure nodes of the FEM mesh,
covering $t \in [0,4]\,\si{s}$. The ramp ends at $t = 1\,\si{s}$, so the shedding is fully
developed over most of the horizon. The control parameters $\theta$ are $U_{\max}$ and the
knot times of the ramp. For training, the grid-based models operate on an internal
$83\times442$ rasterisation of the channel, with the cylinder interior filled by
linear interpolation between the wall nodes and the domain rim by nearest-neighbour extrapolation.

The \ac{kvs} flow requires the full Navier--Stokes equations \eqref{eq:ns}. With
$\nu = 10^{-3}\,\mathrm{m^2/s}$, the Reynolds number
\begin{equation}
  \Reyn = \frac{\bar{U} D}{\nu}, \qquad \bar{U} = \tfrac{2}{3}\,U_{\max},
  \label{eq:re}
\end{equation}
with mean inflow velocity $\bar{U}$, spans $[33, 333]$ across the dataset. Above the critical
value $\Reyn_c \approx 47$, the convective term destabilises the steady flow and leads to
self-sustained periodic vortex shedding, characterised by the Strouhal number
\begin{equation}
  \St = \frac{f_s D}{\bar{U}},
  \label{eq:st}
\end{equation}
with shedding frequency $f_s$. Across the cases above $\Reyn_c$ the shedding period ranges
from $\vKvsShedPeriodMax\,\si{\second}$ just above the threshold to $\vKvsShedPeriodMin\,\si{\second}$ at the highest
inflow, the same order as the one-second ramp. Unlike the film, the wake therefore does not
follow its boundary condition quasi-statically.

The equations are solved fully transiently with an incremental pressure-correction scheme. One
case takes $0.77\,\si{\hour}$ of wall-clock time (Table~\ref{tab:cost}).

\subsection{Splits}
For both datasets, we consider interpolation and extrapolation. For interpolation, the cases are
split randomly with a fixed seed, in a $70/15/15$ ratio of training, validation, and test cases.
This leaves $15$ test cases for \ac{cmp} and $138$ for \ac{kvs}. In the \ac{kvs} case grid, a
test case shares its ramp schedule, at different inflow levels, with training cases. For
extrapolation, \ac{cmp} holds out the $15$ cases with the highest wafer rotation rates, and
\ac{kvs} the $200$ cases with $U_{\max} > 4.0$ ($\Reyn \ge 273$), far above the shedding onset.
The two hold-outs probe different kinds of change by design. The film's Stokes equations are linear
in the boundary velocity, so a higher rotation rate mainly rescales the solution family. The
wake's hold-out instead raises the Reynolds number. There the nonlinearity changes the dynamics
rather than rescaling them.

\section{Surrogate Models}
\label{sec:methods}
We compare a range of prediction models on both datasets. They differ in two respects: the architecture 
used for prediction, and the scheme by which it is deployed in space and time. This section describes the
architectures, then the schemes, and locates every model in Table~\ref{tab:models}.

\subsection{Model architectures}
\label{sec:architectures}

The architectures form a ladder of memory for the control history. A linear map and a
plain \ac{mlp} carry none. Three branch families add that memory in three forms: a gate, a linear state-space recurrence with fixed time scales, and a state-space recurrence
whose time scales follow the input. The three share the deployment scheme and differ only in
the branch, so an answer common to all three belongs to the scheme and not to the
architecture.

\paragraph{Linear baseline (DMDc)}
\Ac{dmdc} fits a single linear map of the latent state by least squares over all training
trajectories jointly \cite{proctor2016dmdc},
\begin{equation}
  \latent(t{+}1) = K\,\latent(t) + L\,c(t), \qquad c(t) = [\,s(t),\, \dot s(t)\,].
  \label{eq:dmdc}
\end{equation}
The fit does not constrain the spectrum of $K$. If eigenvalues fall outside the unit circle,
the rollout diverges, so the fitted operator is constrained to the stable set (spectral radius
$\le 1$). This linear map is the control-aware case of a Koopman linearisation of nonlinear
dynamics \cite{brunton2016koopman,lusch2018koopman}. Parametric variants interpolate such maps
across Reynolds number on the cylinder wake \cite{du2026dmdprom}.

\paragraph{Neural baseline (MLP)}
A \ac{mlp} \cite{rumelhart1986backprop} reads the latent state, the control parameters
$\theta$, the time, and the local ramp value and slope. It is nonlinear but does not learn an
operator.

\paragraph{Operator network (DeepONet)}
A deep operator network (DeepONet) learns the operator directly by factorising the map between
the input and output functions \cite{lu2021deeponet}. A branch network reads the case
inputs, the control parameters $\theta$ and the ramp signal $s(t)$, and returns $p$ coefficients. A trunk
network reads a query $\xi$, a point in space and time for the direct models and the current
latent state with its control values for the autoregressive ones, and
returns $p$ basis values there. Their inner product, plus a bias, is the output at that point,
\begin{equation}
  \mathcal{G}(\theta, s)(\xi) \approx
  \big\langle \mathrm{branch}(\theta, s),\; \mathrm{trunk}(\xi) \big\rangle + b .
  \label{eq:deeponet}
\end{equation}
The trunk is a fixed multilayer perceptron over the query point. We vary only the branch, the
part that processes the case-specific inputs. He et al.\ introduced the
sequential DeepONet, which replaces the naive \ac{mlp} branch with a recurrent neural network
to enhance the model's ability to cope with time-dependent input \cite{he2023sdeeponet}. We
follow this idea and compare three recurrent branches that differ mainly in the form of the
recurrence.

\paragraph{GRU branch}
A recurrent network reads its input sequence step by step. Its hidden state $h_t$ is a vector
that summarises what the network has seen so far. A gated recurrent unit (GRU)
\cite{cho2014gru} reads the case inputs and updates this state at each step as
\begin{equation}
  h_t = z_t \odot h_{t-1} + (1 - z_t) \odot \tilde{h}_t ,
  \label{eq:gru}
\end{equation}
where the update gate $z_t$ and the candidate state $\tilde{h}_t$ are learned nonlinear functions of
$h_{t-1}$ and the input. The gate interpolates between keeping the previous state and
overwriting it. A linear layer maps the final hidden state to the branch coefficients.

\paragraph{State-space branches (S4D and Mamba)}
The two remaining branches replace the gated recurrence with state-space layers. A state-space
layer advances a hidden state along its input sequence $\eta$ with the linear recurrence
\begin{equation}
  h_k = \bar{A}\, h_{k-1} + \bar{B}\, \eta_k, \qquad y_k = 2\,\mathrm{Re}(C\, h_k) + D\,\eta_k.
  \label{eq:ssm}
\end{equation}
This is the discretised form of a continuous-time linear system, with $A$, $B$, $C$, the feedthrough $D$, and the
step size learned. The factor two and the real part follow from the conjugate pairing of the
complex-valued state, of which S4D stores one half. The state recurrence provides a linear alternative for a long memory of
the ramp, without passing information through a gate.
The structured state-space sequence model S4 stacks such layers into a deep network
\cite{gu2022s4}. S4D further simplifies the complicated state matrix carrying the long memory
by making it diagonal and complex-valued, while retaining comparable performance
\cite{gu2022s4d}. A complex-valued state makes each mode of the layer a discretised damped
oscillator~\cite{gu2022s4d}, a natural building block for a periodic wake. Mamba instead makes the state-space matrices functions of the current input \cite{gu2023mamba},
\begin{equation}
  h_k = \bar{A}(\eta_k)\, h_{k-1} + \bar{B}(\eta_k)\, \eta_k.
  \label{eq:mamba}
\end{equation}
It can therefore emphasise or suppress parts of the ramp by their content, and its time
scale can follow the ramp, while S4D applies one fixed recurrence to every ramp. In both cases, the time-averaged output of the last layer,
concatenated with the projected control parameters $\theta$, forms the branch coefficients. We
refer the reader to \cite{gu2023thesis} for the foundations of the model family.

\subsection{Operating schemes}

The models vary in how they represent the state and how they traverse time. Those two dimensions
span Table~\ref{tab:models}. The physics of the two flows fixes them. The film is a linear
response to the momentary control value and carries no state of its own. The wake carries a
phase that is set by its own history, the onset and every shedding period since. The treatment
of time is therefore the axis along which the two flows can separate. The representation is the
second axis. Both flows differ in how far the state can be compressed without losing a
significant share of its variance.

Each cell of the matrix holds one established representative. Within a cell, architectures
differ in capacity and inductive bias but not in what defines the cell. For
the direct cell that is the DeepONet, whose branch reads the control history and whose trunk
reads the query point in space and time. The trunk input follows the regime. On the wake it
reads the normalised position and time, because the field at a point depends on the phase of
the cycle. On the film the field follows the boundary condition at every instant, so time
enters through the branch alone. Its trunk reads the position, the height in the gap, and the
position relative to the wafer centre in polar form. The film is organised around a wafer that
sits at a different offset in every case, and the polar form gives the trunk that geometry
directly. The latent one-shot models use the same branch
architectures as the autoregressive ones, each trained on its own, so that only the treatment of
time varies between the two rows. The latent autoregressive cell decides the regime question and
therefore holds the full ladder of Sec.~\ref{sec:architectures}. The claim that no single
architecture wins both regimes is a claim over the cells this matrix spans, not over every
published architecture.

\paragraph{Representation (direct and latent)}
The models operate either directly on the flow field or on a latent representation of it. A direct
model predicts the field end-to-end. The field DeepONet is direct in this sense. Its basis
is learned jointly with the map rather than fixed beforehand. A latent model reduces the field to a low-dimensional
state, predicts in that space, and decodes back to the field. This keeps the temporal models
small. Encoder and decoder come first, trained once for the wake and computed once in closed
form for the film, and stay frozen. All models of a regime share one encoder per split, so a
difference in performance reflects the architecture, not a co-adapted encoder.

\begin{figure}[t]
\centering
\includegraphics[width=\linewidth]{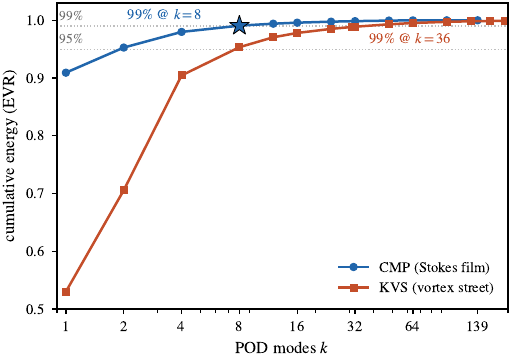}
\caption{Cumulative \ac{pod} energy of the two regimes (train split). \ac{cmp} concentrates 99\%
of the variance in eight modes, the truncation used here. \ac{kvs} reaches 99\%
only at $k=\vKvsPodModesNinetyNine$, with the leading mode at $\vKvsPodModeOnePct\%$, representing a phase-paired limit-cycle structure.}
\label{fig:pod-spectra}
\end{figure}

A \ac{pod} of the training fields gives an idea of how many modes a problem spans
(Fig.~\ref{fig:pod-spectra}). The quasi-static \ac{cmp} film concentrates $99\%$ of its
variance in eight modes. We truncate its \ac{pod} basis at $k=8$, holding $99.1\%$ of the
energy, following standard energy-based truncation~\cite{berkooz1993pod}. Beyond the leading modes the coefficients are
noise-dominated and destabilise the rollout. The oscillatory \ac{kvs} wake spreads its energy across the
coupled mode pairs of a limit cycle and reaches $99\%$ only at $k=\vKvsPodModesNinetyNine$. We therefore compress
it nonlinearly, with a convolutional autoencoder \cite{lee2020convae} of latent dimension 32
reading a single rasterised snapshot with $u$, $v$, and $p$ as its three input channels.
The frozen encoder bounds every latent row. Passing the test fields through encode and decode
alone returns a rollout \ac{rmse} of $\vKvsEncFloorInterpRmse$, a phase-aligned \ac{rmse} of
$\vKvsEncFloorInterpRmsew$, and a shedding power of $\vKvsEncFloorInterpP$ on interpolation. On
extrapolation the rollout \ac{rmse} is $\vKvsEncFloorExtrapRmse$, the phase-aligned \ac{rmse}
$\vKvsEncFloorExtrapRmsew$, and the shedding power $\vKvsEncFloorExtrapP$.
We read every latent row against these values and call them the encoder ceiling for the shedding
power and the encoder floor for the two errors.
Using a different encoder per regime confounds encoder with regime by design, so cross-regime
statements stay at the level of model families, as discussed in Sec.~\ref{sec:limitations}.

\paragraph{Time (autoregressive and one-shot)}
A model either predicts the whole transient in a single pass or advances autoregressively,
feeding its own output back. The autoregressive models start from the encoded initial state
and add a predicted increment to the current state,
\[
  \latent(t{+}1) = \latent(t) + \Delta\latent(t),
\]
except \ac{dmdc}, which advances by the linear map of Eq.~\eqref{eq:dmdc}. This output
feedback is a separate loop from the branch recurrence, which reads the known ramp.

An autoregressive model is trained on its own output, so a small per-step error accumulates
along the rollout. We counter this during training in two ways. The rollout horizon is grown
from a single step to 50 steps, and to 100 for the branch operators on the wake, so the model learns to recover from its own drift.
At evaluation the rollout spans the full transient, beyond the training horizon. A variance
penalty keeps it from settling on the temporal mean, the trivial solution that minimises the
one-step error while erasing the dynamics.

\begin{table}[t]
\centering
\caption{The model set as the design matrix of its two axes, the representation
against the treatment of time. DON abbreviates DeepONet. All models run on both
regimes, and the direct autoregressive cell is left open.}
\label{tab:models}
\footnotesize
\setlength{\tabcolsep}{2.5pt}
\begin{tabular}{@{}l l l@{}}
\toprule
 & one shot & autoregressive \\
\midrule
direct & field DeepONet & --- \\
(full field) & & \\
\addlinespace
latent & S4D one-shot & S4D-DON \\
(frozen encoder) & GRU one-shot & Mamba-DON \\
 & & GRU-DON \\
 & & AR-MLP \\
 & & \ac{dmdc} \\
\bottomrule
\end{tabular}
\end{table}

Table~\ref{tab:models} leaves one combination unfilled, a direct model that advances
autoregressively. We leave this cell open and give the reasons in Sec.~\ref{sec:limitations}.

\subsection{Training and fairness}
\label{sec:training}
Every model minimises a mean squared error on its own representation, the field for the direct
models and the latent state for the latent ones, with the penalty on the temporal variance of the rollout
training as the only addition. No evaluation metric enters any loss, so a lead on a process
target reflects the architecture rather than a directly optimised metric. Every model trains
under a fixed epoch budget, and the best validation checkpoint is reported, chosen on the validation rollout for the
autoregressive models. For the direct operators the validation window covers only the ramp, a
computational shortcut. It does not cause the damping. The last checkpoint of the field DeepONet
damps the wake alike, at $P \le \vKvsFieldDonLastP$ on every interpolation seed, and a mean
squared error over the full horizon would favour the damped prediction further
(Sec.~\ref{sec:metrics}). Each budget was set
heuristically before the benchmark runs, raised until further epochs no longer improved the
best checkpoint, and then held fixed. For the regime champions, some seeds still peaked late in
the budget, so their reported accuracy is a conservative estimate. The learned models are
reported as the mean over three seeds with the standard deviation on the fixed splits of
Sec.~\ref{sec:data}. The deterministic \ac{dmdc} baseline is exact at a single run.

Compute is not equalised across models, but the remaining differences run against our
conclusions. On the wake, the direct field DeepONet trained with the same epoch budget as the
latent operators, at a multiple of their cost per epoch, so its failure to sustain the
oscillation is not an artefact of under-training. Where tuning effort remained unequal, it went
to the models that were behind. The latent one-shot controls received their own search in each
regime and still trail the respective regime champion.

One asymmetry remains. Each champion was developed on its home regime, the direct operator on
the film and the latent pipeline on the wake, and the architectures were reused in the other
regime rather than tuned per cell. Every model is trained from scratch in each regime. What
carries over is the architecture and the loss, with the epoch budget set per regime as above and,
for the three branch operators on the film, the learning rate, the rollout horizon, and the
gradient clip halved for stability. The \ac{mlp} baseline was set separately in each regime. The asymmetry therefore runs both ways. It can shift
margins, but it cannot create the regime flip. The within-wake ablation of
Sec.~\ref{sec:discussion} adds a controlled check, where the branch family and the frozen
encoder are held fixed and only the time treatment changes. All configurations, seeds, and
split definitions are part of the accompanying code.

\section{Evaluation Metrics}
\label{sec:metrics}
\subsection{Shared aspects}
\label{sec:mirrored-metrics}
The simplest summary of a surrogate's accuracy is a single pointwise \ac{rmse} against the reference. For the
wake this is the rollout \ac{rmse}, the pointwise error over the $4\,\si{\second}$ trajectory, with $N$ evaluation nodes
$\xvec$, $T$ frames $t$, the velocity components $i$, and the $C$ test cases $c$,
\begin{equation}
  \mathrm{RMSE} = \sqrt{\tfrac{1}{2NTC}\textstyle\sum_{c,\xvec,t,i}\big(\hat u_i(\xvec,t) - u_i(\xvec,t)\big)^2}, \label{eq:rollout}
\end{equation}
in raw velocity units, since the two wake components share a common scale. For the film the
aggregate is evaluated over the two in-plane velocity components at every node of the full
three-dimensional domain, each normalised by its global standard deviation $\sigma_i$,
\begin{equation}
  \mathrm{nRMSE}_{\mathrm{ip}} = \sqrt{\tfrac{1}{2NTC}\textstyle\sum_{c,\xvec,t,i\in\{x,y\}}\big((\hat u_i(\xvec,t) - u_i(\xvec,t))/\sigma_i\big)^2}. \label{eq:cmp_rmse}
\end{equation}
The through-gap component is excluded from the aggregate. It lies more than six orders of magnitude
below the in-plane flow. Under per-component normalisation it would therefore enter the
aggregate as a noise floor with the same weight as the two physical components. We report it descriptively with the results instead.

For a transient flow the aggregate misleads, both as a training objective and as a
score. As an objective it rewards damping, because a smooth, low-amplitude
guess is rarely far from the truth~\cite{lippe2023pderefiner}. As a score it rates a wrong structure, a damped amplitude, and invented motion all as the
same mild deviation~\cite{mohan2024illusion}, while a prediction that is merely displaced in
time can score worse than the damped one. It can even rank
predictions in the opposite order to their physical fidelity when several outcomes are
admissible~\cite{baattrup2026}.

The reason is elementary. Take a probe trace of the developed wake, $v(t) = A\sin(\omega t)$.
A prediction that damps the oscillation to its time mean, $\hat v = 0$, has a pointwise
\ac{rmse} of $A/\sqrt{2}$. A prediction with the correct amplitude, frequency, and onset,
displaced by a phase $\varphi$, has $\sqrt{2}\,A\,\lvert\sin(\varphi/2)\rvert$. The damped
prediction wins as soon as $\varphi$ exceeds a sixth of a period. At half a period the perfect
oscillation scores exactly twice as badly as the flat line. That case is the wake's mirror
branch, taken up in Sec.~\ref{sec:z2}. The
pointwise norm penalises a displaced oscillation twice, once where it predicts motion the
reference does not show at that instant and once where it misses the motion the reference
shows. Forecast verification knows this as the double
penalty~\cite{ebert2008fuzzy,gilleland2009intercomparison}, and it is the reason a damped
surrogate can outrank a physically correct one. The metric suite below is built against this
defect.

\begin{figure*}[t]
\centering
\includegraphics[width=\textwidth]{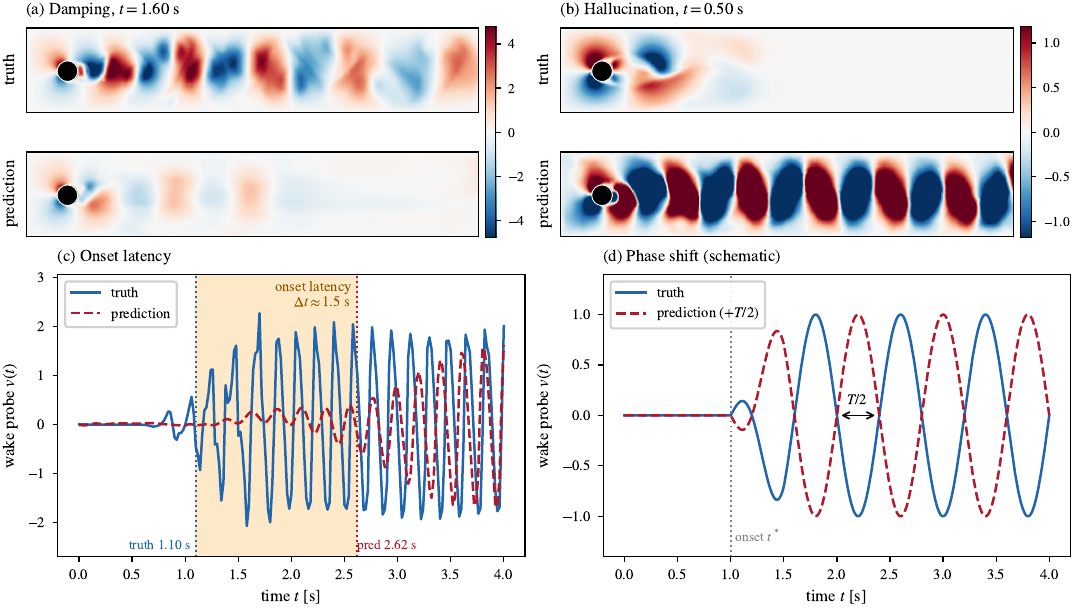}
\caption{The four \ac{kvs} failure modes that a pointwise \ac{rmse} folds together.
\textbf{(a)}~Damping. The prediction washes the developed wake out toward zero.
\textbf{(b)}~Hallucination. The rollout invents a full vortex street at $t=0.50$\,s, at
several times the amplitude the flow carries at that instant and well before the true
street is established.
\textbf{(c)}~Onset latency. The prediction stays smooth while the
truth already sheds, a delayed onset. The marks are the detector onsets, the
times shedding is established. \textbf{(d)}~Phase shift (schematic). The prediction oscillates at the correct frequency,
amplitude, and onset, displaced by half a period, the largest possible phase error on the limit
cycle. Top
row (amplitude) shows $v$-field stills, truth over prediction. Bottom row (timing) shows
wake-probe $v(t)$ traces. Panels (a) to (c) are single-seed example predictions from the
benchmarked model set, chosen to isolate each mode.}
\label{fig:kvs-failures}
\end{figure*}

Fig.~\ref{fig:kvs-failures} shows four wake failure modes that the pointwise norm folds
together, on example cases.
The first three, damping, hallucination, and onset latency, are scored below. The fourth, a phase
shift (d), displaces the oscillation in time while leaving frequency, amplitude, and onset
unchanged. It is not scored as a failure mode of its own. The pointwise norm penalises it, as
derived above, and the phase-aligned \ac{rmse} removes it from the structural score.

We evaluate each surrogate model against a small, fixed set of metrics. Each metric isolates one
component that the pointwise norm folds into a single number. Five core aspects of a prediction
organise the suite. The first three are posed in both regimes, as the same physical question
answered with a regime-appropriate instrument, and the last two exist only in the wake:
\begin{enumerate}
  \item the \emph{field} as a whole, its overall accuracy,
  \item the \emph{spatial structure}, whether the pattern is correct,
  \item any \emph{invented motion} the reference does not contain,
  \item the \emph{fluctuation amplitude}, whether the model damps the dynamics, and
  \item the \emph{timing} of events.
\end{enumerate}
The aspects follow from the defect of the
pointwise norm, sharpened by early experiments and by visual inspection of the predicted fields.
Deriving the metrics from the physics rather than from the results keeps
them from being chosen, even unintentionally, to favour a particular model~\cite{gelman2014forking}.
We report seven metrics for the wake and five for the film. Throughout, $\mu_{x}[\cdot]$,
$\sigma_{x}[\cdot]$, and $\sigma^2_{x}[\cdot]$ denote the mean, standard deviation, and variance
over the subscripted variable.

For the wake, the field-level measures, like the rollout \ac{rmse} and the hallucination fraction,
are evaluated mesh-natively on a fixed subset of $8{,}192$ velocity nodes frozen across all
models. The subset was drawn once, uniformly over the mesh nodes, so its density follows the
mesh grading. The oscillatory metrics,
frequency, phase, and amplitude, read the prediction $\hat v$
against the truth $v$ at twelve fixed probes along the channel.

For the film, the remaining metrics are evaluated on the inner wafer, the central region
$r \le 0.9\,r_{\mathrm{wafer}}$, denoted $\Omega_{\mathrm{in}}$ in the equations below. The rim is excluded for a numerical reason. Past the wafer edge
the domain continues as a traction-free annulus, but the prescribed rotation ends abruptly at
the edge, and the steep gradient there is the least certain part of the reference. Metrics
reading the rim would score a surrogate against reference artefacts rather than against the flow.

\paragraph{Spatial structure}
For the wake a pure phase lag must not count as a structural error, so the phase-aligned \ac{rmse}
aligns prediction and truth in time by \ac{dtw} before scoring,
\begin{equation}
  \rmsew = \mu_{c,p}\Big[ \min_{\gamma \in \mathrm{DTW}} \sqrt{\tfrac{1}{|\gamma|}\textstyle\sum_{(i,j)\in\gamma}
    \big(\hat v_{c,p}(t_i) - v_{c,p}(t_j)\big)^2} \Big], \label{eq:rmsew}
\end{equation}
with $\gamma$ a warp path, a set of index pairs $(i,j)$ that matches predicted to true time samples,
and $|\gamma|$ its length. The path is monotone and anchored at both ends, which removes a phase
lag without aligning away genuinely distinct structures. Each trace is aligned by its own warp
path. The path is found by dynamic programming on the absolute deviation between the two
traces, with the symmetric step pattern that charges the diagonal step twice. The score
then reads the squared deviations along that path, so Eq.~\eqref{eq:rmsew} states the reported
quantity rather than the cost the search minimises. The phase-aligned \ac{rmse} reads the twelve wake probes rather than the field nodes, so its
values are not directly comparable to the rollout \ac{rmse}.

For the film, where nothing oscillates, the structural question is whether the predicted
\ac{cumwss} accumulates in the right places on the wafer. The Pearson correlation
$\mathrm{corr}(\hat{\mathcal{S}},\mathcal{S})$ over the inner-wafer nodes scores this. A high
correlation means the high- and low-shear zones sit where the reference puts them. It is blind
to magnitude.

\paragraph{Invented motion}
For the wake the true flow carries no shedding before onset (Fig.~\ref{fig:kvs-failures}b), so any
predicted motion there is classified as hallucination, quantified by two measures. The first, $H$,
measures the amount of invented motion. Deciding when a velocity value counts as motion takes a
threshold. We set it to a tenth of the peak steady velocity of the true flow over all test
cases $c$,
\begin{equation}
  \tau = 0.1 \cdot \max\nolimits_{c,\, t \ge 2.5\,\si{\second},\, \xvec,\, i} |u_i(\xvec,t)|, \label{eq:halluc_tau}
\end{equation}
one value shared across the test set, so that $H$ is comparable across cases and models. For each
case $c$ we then collect every velocity value in the fixed early window $t < 0.5\,\si{\second}$,
one sample per frame, evaluation node, and component, into the set $\mathcal{P}_c$. A sample is hallucinated when
the truth rests below the threshold while the prediction moves above it. $H$ is the fraction of
hallucinated samples, averaged over the test cases,
\begin{equation}
  H = \mu_c\Big[ \tfrac{1}{|\mathcal{P}_c|}
    \textstyle\sum_{(t,\xvec,i)\in\mathcal{P}_c}
    \mathbf{1}\big[\,|u_i(\xvec,t)| < \tau \,\wedge\, |\hat u_i(\xvec,t)| > \tau\,\big] \Big]. \label{eq:halluc_frac}
\end{equation} The second
is whether that motion resembles a K\'arm\'an street or is unstructured noise. $\alpha_{\mathrm{rms}}$ projects
the pre-onset residual onto the developed-wake template
$\Phi = \uvec(t_{\mathrm{end}}) - \mu_{t \ge 2.5\,\si{\second}}[\uvec(t)]$,
\begin{equation}
  \alpha(t) = \frac{\big\langle\, \hat\uvec(t) - \mu_{t'\in\mathrm{pre}}[\uvec(t')],\ \Phi \,\big\rangle}
                   {\lVert \Phi \rVert^2},
  \qquad
  \alpha_{\mathrm{rms}} = \sqrt{\mu_{t < 0.5\,\si{\second}}\big[ \alpha(t)^2 \big]},
    \label{eq:alpharms}
\end{equation}
with $\langle\cdot,\cdot\rangle$ the inner product over both velocity components on an
$83\times442$ raster, so
$\alpha = k$ means $k$ times the settled vortex pattern is already present. The raster
spans the full channel of Sec.~\ref{sec:data}. It is filled from the evaluation nodes by linear
interpolation on their Delaunay triangulation, and by the nearest node where that interpolation
is undefined, the same raster for prediction and truth.
Read together they separate the cases. High $H$ with low $\alpha_{\mathrm{rms}}$ is incoherent noise.
Low $H$ with high $\alpha_{\mathrm{rms}}$ is a faint but coherent pseudo-wake.

For the film the same failure is a single
metric. Once the ramp saturates, the Stokes field stops evolving in time. Any temporal variation a
model predicts in the developed phase therefore has no counterpart in the true flow. Because the true
developed-phase variance is near zero, an absolute ratio is ill-conditioned. The spurious
unsteadiness is instead measured as the excess over truth on the two-plane shear speed
$w(\xvec,t) = \lvert \uvec_{\text{top}}(\xvec,t) - \uvec_{\text{bot}}(\xvec,t) \rvert$,
\begin{equation}
  D = \frac{\mu_{\xvec\in\Omega_{\mathrm{in}}}\big[\sigma_{t \ge 1.5\,\si{\second}}[\hat w(\xvec,t)]\big]}{\mu_{\xvec\in\Omega_{\mathrm{in}}}\big[\mu_{t \ge 1.5\,\si{\second}}[\hat w(\xvec,t)]\big]}
    - \frac{\mu_{\xvec\in\Omega_{\mathrm{in}}}\big[\sigma_{t \ge 1.5\,\si{\second}}[w(\xvec,t)]\big]}{\mu_{\xvec\in\Omega_{\mathrm{in}}}\big[\mu_{t \ge 1.5\,\si{\second}}[w(\xvec,t)]\big]}.
    \label{eq:cmp_drift}
\end{equation}
Even the settled true field retains a small residual variation, so the truth term is subtracted and
$D$ reports only the excess a model adds over the true flow.

\subsection{Regime-exclusive aspects}
\label{sec:exclusive-metrics}
\paragraph{Process target}
The film's exclusive aspect is its process target, the \ac{cumwss}, the time-integrated
wall shear stress. It is scored per case by its relative error over the inner wafer,
\begin{equation}
  \mathrm{cwL_2} = \frac{\lVert \hat{\mathcal{S}} - \mathcal{S} \rVert_2}{\lVert \mathcal{S} \rVert_2}. \label{eq:cmp_cumwss}
\end{equation}
It also carries the magnitude that $\mathrm{corr}(\hat{\mathcal{S}},\mathcal{S})$ deliberately
ignores. The bulk field error and the \ac{cumwss} need not rank the models by the same
margin, a contrast we develop in Sec.~\ref{sec:cmp}.

\paragraph{Fluctuation amplitude}
A model may relax the wake toward its time-mean (Fig.~\ref{fig:kvs-failures}a). The shedding power
$P$ measures how much of the true fluctuation energy the prediction retains. This energy is captured by the temporal variance of a probe trace, so
\begin{equation}
  P = \frac{\sum_{c,p} \sigma^2_{t \ge 2\,\si{\second}}[\hat v_{c,p}(t)]}{\sum_{c,p} \sigma^2_{t \ge 2\,\si{\second}}[v_{c,p}(t)]}, \label{eq:varratio}
\end{equation}
pooled over the test cases $c$ and the twelve probes $p$, which weights each probe by its true
signal. It falls below one when the wake is damped.

\paragraph{Onset timing}
\label{sec:onset}
Onset is the event in which the wake breaks its top-bottom symmetry
(Fig.~\ref{fig:kvs-failures}c). Over three stripe
pairs straddling the cylinder axis we average $v$ over the
evaluation nodes in the sets $\mathrm{top}_k$ above and
$\mathrm{bot}_k$ below the axis. Each stripe is $0.04$ wide and $0.10$ high, and the three
stations sit $0.02$ upstream and $0.02$ and $0.06$ downstream of the cylinder centre. While the wake is symmetric, the two means cancel
($v_{\text{top}} = -v_{\text{bot}}$), so the stripe signal
\begin{equation}
  s_k(t) = \mu_{\xvec\in\mathrm{top}_k}[v(\xvec,t)] + \mu_{\xvec\in\mathrm{bot}_k}[v(\xvec,t)]
  \label{eq:ysym}
\end{equation}
rests at zero until symmetry breaks. Shedding makes $s_k$ oscillate, so the detector reads its
frame-to-frame fluctuation, the absolute difference between consecutive frames,
summarised by a centred rolling median over $13$ frames, about one shedding period. Onset
$t^*$ is the first time this fluctuation exceeds a tenth of its developed level on
the truth, the same bar for prediction and truth. The developed level is the $95$th
percentile of the fluctuation over the trajectory. The bar counts as exceeded when it is passed
in at least two of the three stripes for three consecutive frames, within $0.4$ to
$3.5\,\si{\second}$. The onset latency
$\Delta t^{*} = |t^*_{\hat v} - t^*_v|$ is the timing error between prediction and truth.
The detector parameters shift the absolute onset times but affect the model comparison only
marginally. The interval from first visible motion to established shedding is itself short and
nearly model-independent, so $\Delta t^{*}$ measures a start offset of the instability rather
than a build-up penalty. $\Delta t^{*}$ is defined only where the detector finds an onset
on both the truth and the prediction. For a collapsed
model that never sheds there is no onset to hit or miss, and it is reported as undefined.

$\Delta t^{*}$ also resolves the ambiguity the structural score $\alpha_{\mathrm{rms}}$ leaves
open. A high $\alpha_{\mathrm{rms}}$ on its own cannot differentiate an early-but-real shedder
from a hallucinator. A small $\Delta t^{*}$ means the model sheds correctly, only a little
early, while structure with no matching onset is genuine hallucination. Sec.~\ref{sec:kvs}
reports both.

\subsection{Domain indicators}
\label{sec:kpi-metrics}
Each regime has its own domain indicator, with no cross-regime counterpart. For the wake it is the
Strouhal error,
\begin{equation}
  \fftfreq = \mu_{c,p}\Big[ \tfrac{|\hat f_p - f_p|}{f_p} \Big], \label{eq:fft}
\end{equation}
with $f_p$ the dominant shedding frequency of $v$ at probe $p$ over the developed phase
$t \ge 2\,\si{\second}$. The frequency is read from the spectrum of the mean-removed trace,
Hann-windowed and zero-padded fourfold, with a parabolic refinement of the peak through the log
power of the peak bin and its two neighbours. The resolution therefore lies well below the
$0.5\,\si{\hertz}$ bin of the raw window. The mean runs over every case-probe pair whose true trace has a
spectral peak away from zero frequency. This condition drops $\vKvsFftDroppedInterp$ of the
$\vKvsFftPairsInterp$ pairs on the interpolation split and none on extrapolation. Pairs at
quiet probes of sub-critical cases stay in and add a small floor shared by all models. A wrong
frequency cannot be repaired by any amplitude or phase correction. For the film it is the \ac{wiwnu}, the spatial coefficient of variation of the \ac{cumwss},
reported as the difference between prediction and truth,
\begin{equation}
  \Delta\mathrm{WIWNU} = \frac{\sigma_{\xvec\in\Omega_{\mathrm{in}}}[\hat{\mathcal{S}}(\xvec)]}{\mu_{\xvec\in\Omega_{\mathrm{in}}}[\hat{\mathcal{S}}(\xvec)]}
    - \frac{\sigma_{\xvec\in\Omega_{\mathrm{in}}}[\mathcal{S}(\xvec)]}{\mu_{\xvec\in\Omega_{\mathrm{in}}}[\mathcal{S}(\xvec)]}.
    \label{eq:cmp_wiwnu}
\end{equation}

The rollout \ac{rmse}, $\mathrm{nRMSE}_{\mathrm{ip}}$, the Strouhal error, and $P$
pool all cases before aggregation, and the phase-aligned \ac{rmse} averages over all case-probe pairs. For the Strouhal error and $P$ the pooling is also what keeps
them well defined when sub-critical test cases do not shed and per-case frequencies or ratios
are ill-posed. All remaining metrics are means over test cases. Every
metric is reported as the mean over three seeds with its standard deviation, on the fixed splits
of Sec.~\ref{sec:data}.

\subsection{Computational cost and amortisation}
\label{sec:cost}

The motivation for a surrogate is speed, so we report cost alongside accuracy. The raw
per-query speedup is the ratio of solver to surrogate wall-clock,
\begin{equation}
  S = \frac{t_{\mathrm{FEM}}}{t_{\mathrm{infer}}}, \label{eq:speedup}
\end{equation}
where $t_{\mathrm{FEM}}$ is the finite-element wall-clock for one case and
$t_{\mathrm{infer}}$ is the warm per-case inference latency, both measured to the same usable field output.

This ratio alone overstates the practical gain. A data-driven surrogate first pays an offline cost, the
$N_{\mathrm{train}}$ reference simulations that form the training set together with the training itself.
For a latent pipeline the training term $t_{\mathrm{train}}$ covers both stages of
Sec.~\ref{sec:training}, the encoder fit and the dynamics fit. The search that selected each
configuration is not charged. The offline cost therefore prices the reproduction of the final
pipeline, not its development.
For a campaign of $Q$ queries the two totals are
\begin{equation}
\begin{aligned}
  C_{\mathrm{FEM}}  &= Q\,t_{\mathrm{FEM}}, \\
  C_{\mathrm{surr}} &= \underbrace{N_{\mathrm{train}}\,t_{\mathrm{FEM}} + t_{\mathrm{train}}}_{\mathrm{offline}}
                       + Q\,t_{\mathrm{infer}},
\end{aligned}
\label{eq:amort}
\end{equation}
which cross at the break-even query count
\begin{equation}
  Q^{*} = \frac{N_{\mathrm{train}}\,t_{\mathrm{FEM}} + t_{\mathrm{train}}}{t_{\mathrm{FEM}} - t_{\mathrm{infer}}}.
  \label{eq:breakeven}
\end{equation}
In the limit $t_{\mathrm{infer}} \ll t_{\mathrm{FEM}}$, the break-even count reduces to
$Q^{*} \approx N_{\mathrm{train}} + t_{\mathrm{train}}/t_{\mathrm{FEM}}$. We report $Q^{*}$ next to $S$,
following the reporting recommendation of McGreivy and Hakim~\cite{mcgreivy2024weak},
since the break-even count, not the raw ratio, decides whether a surrogate pays off for a given workload.

The two wall-clocks come from different hardware, the solver on CPU and the surrogate on GPU, so
$S$ compares deployments rather than devices. We report both wall-clocks alongside the ratio, and
take the solver at its production core count rather than a single core, so the reference is not
artificially slow.

\section{Results}
\label{sec:results}

Every model of Sec.~\ref{sec:methods} is scored through the common evaluation pipeline of
Sec.~\ref{sec:metrics}, on the interpolation and extrapolation splits of Sec.~\ref{sec:data}. Rows in the
tables are grouped by their design-matrix cell of
Sec.~\ref{sec:methods}, representation against time treatment. We present the two regimes in turn
and discuss them together in Sec.~\ref{sec:discussion}.

\subsection{CMP (film)}
\label{sec:results-cmp}
\label{sec:cmp}

\begin{figure}[t]
\centering
\includegraphics[width=\linewidth]{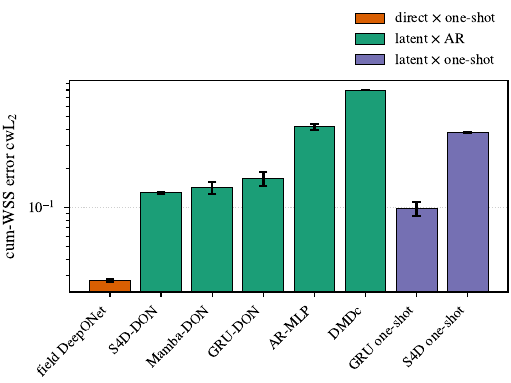}
\caption{\ac{cmp} model comparison on the interpolation split, on the process target
$\mathrm{cwL}_2$ (log scale, lower is better). The direct one-shot operator leads the closest
latent variant by more than a factor of three. Error bars show the standard deviation over
three seeds where the model is stochastic.}
\label{fig:cmp-champions}
\end{figure}

Table~\ref{tab:cmp-results} summarises the five-metric suite, and
Fig.~\ref{fig:cmp-champions} compares the models on the process target.

\begin{table*}[t]
\centering
\caption{\ac{cmp} (film) results, the mean over three seeds with the seed standard deviation as a
subscript. Columns are the in-plane bulk field error
$\mathrm{nRMSE}_{\mathrm{ip}}$ (Eq.~\eqref{eq:cmp_rmse}), the spatial correlation $\mathrm{corr}(\hat{\mathcal{S}},\mathcal{S})$, the excess unsteadiness $D$, the process target
$\mathrm{cwL}_2$ (relative \ac{cumwss} error, Eq.~\eqref{eq:cmp_cumwss}), and the within-wafer non-uniformity
error $\Delta\mathrm{WIWNU}$. Lower is better except the correlation (higher) and $\Delta\mathrm{WIWNU}$ (closer
to zero). Bold marks the best \ac{cumwss}. Rows are grouped by design-matrix cell, not by rank.
$^\dagger$\ac{dmdc} is deterministic, reported once.}
\label{tab:cmp-results}
\footnotesize
\begin{tabular}{@{}ll ccccc@{}}
\toprule
& Model & $\mathrm{nRMSE}_{\mathrm{ip}}$ & $\mathrm{corr}(\hat{\mathcal{S}},\mathcal{S})$ & $D$ & $\mathrm{cwL}_2$ & $\Delta\mathrm{WIWNU}$ \\
\midrule
\multicolumn{7}{@{}l}{\emph{Interpolation}}\\
direct $\times$ one-shot & field DeepONet        & \rowCmpInterpFieldDON \\
latent $\times$ AR       & S4D-DON                & \rowCmpInterpSFourD \\
                         & Mamba-DON              & \rowCmpInterpMamba \\
                         & GRU-DON                & \rowCmpInterpGRU \\
                         & AR-MLP                 & \rowCmpInterpARMLP \\
                         & \ac{dmdc}$^\dagger$         & \rowCmpInterpDMDc \\
latent $\times$ one-shot & GRU one-shot           & \rowCmpInterpOneshotGRU \\
                         & S4D one-shot           & \rowCmpInterpOneshotSFourD \\
\addlinespace
\multicolumn{7}{@{}l}{\emph{Extrapolation}}\\
direct $\times$ one-shot & field DeepONet        & \rowCmpExtrapFieldDON \\
latent $\times$ AR       & S4D-DON                & \rowCmpExtrapSFourD \\
                         & Mamba-DON              & \rowCmpExtrapMamba \\
                         & GRU-DON                & \rowCmpExtrapGRU \\
                         & AR-MLP                 & \rowCmpExtrapARMLP \\
                         & \ac{dmdc}$^\dagger$         & \rowCmpExtrapDMDc \\
latent $\times$ one-shot & GRU one-shot           & \rowCmpExtrapOneshotGRU \\
                         & S4D one-shot           & \rowCmpExtrapOneshotSFourD \\
\bottomrule
\end{tabular}
\end{table*}

On interpolation the in-plane bulk error ranks the direct operator first at
$\mathrm{nRMSE}_{\mathrm{ip}} = 0.23$, ahead of the latent variants from $0.33$ to $0.54$. The
two baselines separate, at $1.0$ for the \ac{mlp} and $2.4$ for \ac{dmdc}. The excluded through-gap
component sits at a normalised error near one for every model, the level of predicting its
mean. Under extrapolation it rises above one, worse than the mean. The spatial
correlation lies at $0.999$ for the direct operator and $0.90$ to
$0.93$ for the recurrent models. The excess unsteadiness $D$ stays near zero for the direct and
recurrent models and rises to $0.078$ for the \ac{mlp} and $0.575$ for the linear baseline. On
the process target the direct operator leads at $\mathrm{cwL}_2 = 0.027$, the GRU one-shot
follows at $0.099$, the recurrent latent-AR models sit at $0.13$ to $0.17$, and the two
baselines fall off at $0.42$ and $0.80$. $\Delta\mathrm{WIWNU}$ stays below $0.06$ in magnitude
for every model except the S4D one-shot and the linear baseline.

Two readings deserve attention. The bulk error and the process target agree on the leading
model but not on its margin, the direct operator ahead of the GRU-DON by a factor of $\vFilmIpFactor$ on
$\mathrm{nRMSE}_{\mathrm{ip}}$ and by a factor of six on $\mathrm{cwL}_2$. And the GRU one-shot, which drops the
autoregressive feedback, beats every latent-AR variant on the process target.

\begin{figure*}[t]
\centering
\includegraphics[width=\textwidth]{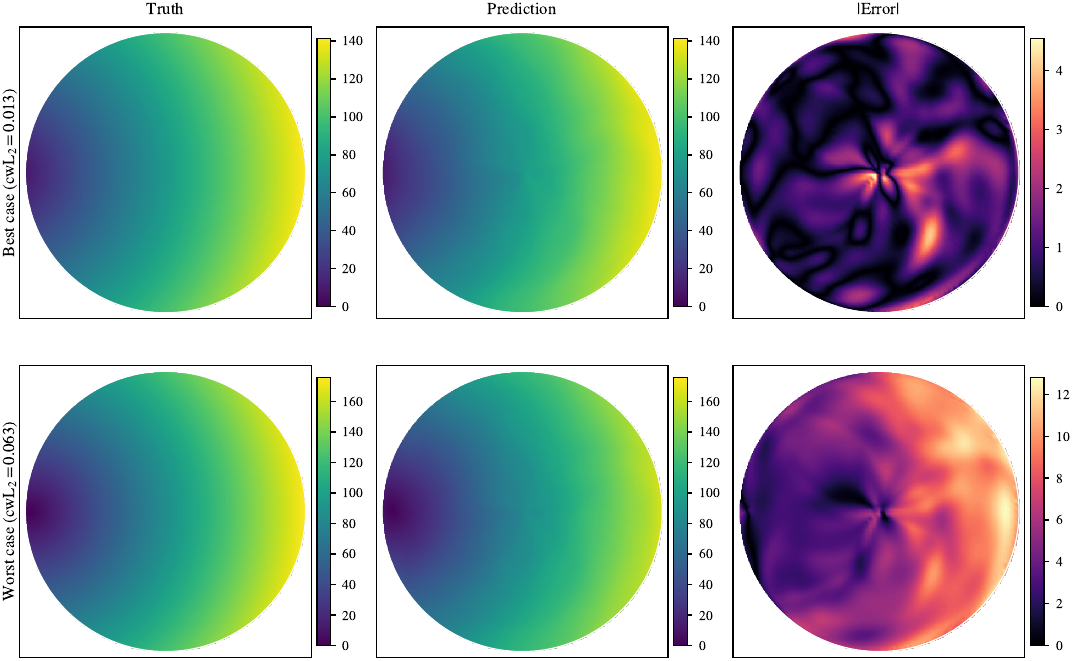}
\caption{\ac{cmp} \ac{cumwss} field of the direct operator (field DeepONet) against truth on its best
and worst interpolation case by $\mathrm{cwL}_2$. Maps show the inner wafer ($r \le 0.9\,r_{\mathrm{wafer}}$,
Sec.~\ref{sec:metrics}), and the truth and prediction columns share one scale per row. The error
column carries its own scale per row. On the best case
the operator reproduces the footprint almost exactly. On the worst case the deviation stays
smooth and large-scale. It is largest along the side of the wafer farthest from the pad's
rotation axis. There the sliding speed
and with it the accumulated stress are highest, so the deviation follows the stress itself rather
than any localised structure.}
\label{fig:cmp-wss}
\end{figure*}

Extrapolation preserves both readings. The direct operator holds the process target at $0.032$,
the GRU one-shot follows at $0.091$, and the recurrent models sit at $0.12$ to $0.15$, while on
$\mathrm{nRMSE}_{\mathrm{ip}}$ the direct operator leads as well, at $0.31$ against $0.37$
to $0.69$ for the recurrent and one-shot latent variants. The
spatial correlation stays at $0.92$ or above for the direct and
recurrent models. The \ac{mlp} baseline breaks down, with $\mathrm{cwL}_2$ rising to $1.77$ and
$\mathrm{nRMSE}_{\mathrm{ip}}$ to $37.7$ in the seed mean. Fig.~\ref{fig:cmp-wss} shows the \ac{cumwss}
field of the direct operator against truth on its best and worst interpolation case.

\begin{figure}[t]
\centering
\includegraphics[width=\linewidth]{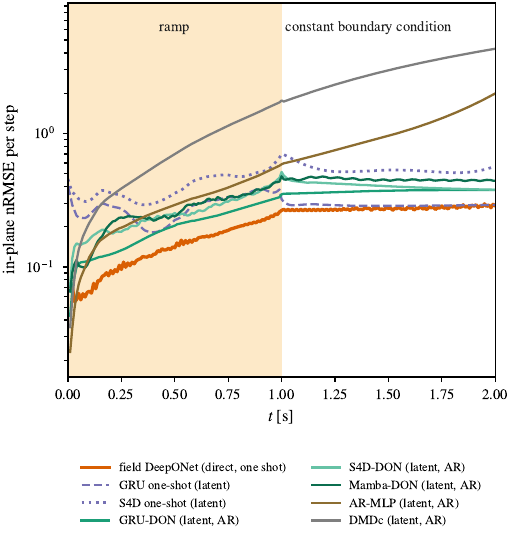}
\caption{Film error along the trajectory on the interpolation split. In-plane
$\mathrm{nRMSE}$ per time step, normalised as in Table~\ref{tab:cmp-results}, mean over the test cases and
the three seeds, from each model's first predicted step. The autoregressive rows start from the
encoded true initial state. The shaded interval is the ramp, after which the boundary condition
is constant.}
\label{fig:cmp-rollout-growth}
\end{figure}

Fig.~\ref{fig:cmp-rollout-growth} reads the film error along the trajectory. The direct
operator's error rises through the ramp in step with the growing field and stays flat once
the boundary condition is constant, at $\vCmpGrowDirectRampEnd$ at the end of the ramp and
$\vCmpGrowDirectEnd$ at the end. Every autoregressive row starts at the encoder floor of
$\vCmpGrowArFloor$, well below the direct operator, and crosses it during the ramp. By the end
of the ramp the best of them, the GRU branch, stands at $\vCmpGrowGruRampEnd$. In the constant
phase the two rows without a memory of the ramp keep growing, the \ac{mlp} by a
factor of $\vCmpGrowMlpRatio$ and the linear map by $\vCmpGrowDmdcRatio$, while the three
recurrent branches plateau at their ramp-end error and never return to the level of the
direct operator.

\subsection{KVS (wake)}
\label{sec:results-kvs}
\label{sec:kvs}

\begin{figure*}[t]
\centering
\includegraphics[width=\textwidth]{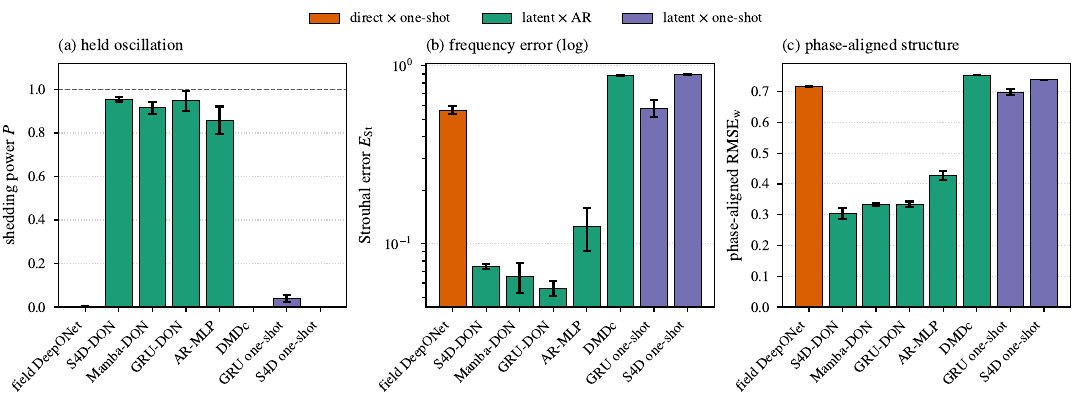}
\caption{\ac{kvs} model comparison on the interpolation split. (a)~Shedding power $P$, closer to one
is better. (b)~Strouhal frequency error $\fftfreq$ (log scale). (c)~Phase-aligned \ac{rmse}
$\rmsew$. The trained latent-AR models separate from the damping models on all three axes. Error bars
show the standard deviation over three seeds where the model is stochastic.}
\label{fig:kvs-champions}
\end{figure*}

Table~\ref{tab:kvs-results} summarises the seven-metric suite for \ac{kvs}, and
Fig.~\ref{fig:kvs-champions} compares the models on the held shedding power, the frequency
error, and the phase-aligned structure.

\begin{table*}[t]
\centering
\caption{\ac{kvs} (wake) results, the mean over three seeds with the seed standard deviation as a subscript.
Columns are the rollout \ac{rmse}, the phase-aligned
$\rmsew$, the pre-onset hallucination fraction $H$, the pseudo-wake amplitude $\alpha_{\mathrm{rms}}$,
the shedding power $P$, the onset latency $\Delta t^{*}$ in seconds, and the Strouhal error $\fftfreq$.
Lower is better except $P$ (closer to one). Bold marks the leading shedding power. Rows are grouped by
design-matrix cell, not by rank. $^\dagger$\ac{dmdc} is deterministic, reported once. A dash marks an onset
undefined for a model whose prediction stays below the onset threshold on more than half of the test cases.}
\label{tab:kvs-results}
\footnotesize
\setlength{\tabcolsep}{2.2pt}
\begin{tabular}{@{}ll ccccccc@{}}
\toprule
& Model & \ac{rmse} & $\rmsew$ & $H$ & $\alpha_{\mathrm{rms}}$ & $P$ & $\Delta t^{*}$ & $\fftfreq$ \\
\midrule
\multicolumn{9}{@{}l}{\emph{Interpolation}}\\
direct $\times$ one-shot & field DeepONet         & \rowKvsInterpFieldDON \\
latent $\times$ AR       & S4D-DON                & \rowKvsInterpSFourD \\
                         & Mamba-DON              & \rowKvsInterpMamba \\
                         & GRU-DON                & \rowKvsInterpGRU \\
                         & AR-MLP                 & \rowKvsInterpARMLP \\
                         & \ac{dmdc}$^\dagger$         & \rowKvsInterpDMDc \\
latent $\times$ one-shot & GRU one-shot           & \rowKvsInterpOneshotGRU \\
                         & S4D one-shot           & \rowKvsInterpOneshotSFourD \\
\addlinespace
\multicolumn{9}{@{}l}{\emph{Extrapolation}}\\
direct $\times$ one-shot & field DeepONet         & \rowKvsExtrapFieldDON \\
latent $\times$ AR       & S4D-DON                & \rowKvsExtrapSFourD \\
                         & Mamba-DON              & \rowKvsExtrapMamba \\
                         & GRU-DON                & \rowKvsExtrapGRU \\
                         & AR-MLP                 & \rowKvsExtrapARMLP \\
                         & \ac{dmdc}$^\dagger$         & \rowKvsExtrapDMDc \\
latent $\times$ one-shot & GRU one-shot           & \rowKvsExtrapOneshotGRU \\
                         & S4D one-shot           & \rowKvsExtrapOneshotSFourD \\
\bottomrule
\end{tabular}
\end{table*}

On interpolation the rollout \ac{rmse} places the S4D and Mamba branches first, at
$0.34$ and $0.33$, followed by the GRU branch at $0.44$. The field DeepONet follows at $0.66$,
ahead of the latent one-shots at $0.70$ and $0.72$, the \ac{mlp} baseline at $0.94$, and the
linear baseline at $1.38$. The
phase-aligned $\rmsew$ keeps the trained latent-AR models ahead, $0.26$ to $0.44$,
against $0.72$ to $0.77$ for the rest. Pre-onset hallucination spans $H = 0.003$ to $0.022$
across the trained latent-AR models. The GRU branch stays at $0.003$, on the same level as the
field DeepONet, while the S4D branch moves most in the seed mean at $0.022$ and the linear baseline
reaches $0.042$.
The pseudo-wake amplitude stays at or below $\alpha_{\mathrm{rms}} = 0.021$ for every trained model. The shedding power separates the suite
most sharply. The latent-AR models hold the oscillation, the Mamba, S4D, and GRU branches at
$P = 0.90$ to $0.91$ and the \ac{mlp} at $0.84$, while the direct, linear, and latent
one-shot models fall to $0.03$ and below. The three recurrent branches reach the encoder ceiling of
$P = \vKvsEncFloorInterpP$ (Sec.~\ref{sec:methods}).
Among the shedders the onset latency is $0.09$ to $0.10\,\si{\second}$ for the three recurrent branches and
$0.39\,\si{\second}$ for the \ac{mlp}. For the field DeepONet and the GRU one-shot the detector returns $0.12$ and
$0.15\,\si{\second}$. Both reproduce the initial transient and then damp it away. These values therefore report
when that transient started, not a sustained onset. The Strouhal error follows the shedding split, $0.04$
to $0.17$ for the latent-AR models against $0.56$ to $0.89$ for the rest.

Three rankings deserve attention. Under extrapolation the rollout \ac{rmse} and the shedding
power sort the models in
opposite directions, the damped field DeepONet ahead on the aggregate, $1.10$ against $1.29$,
while holding none of the oscillation, $P = 0.002$ against $0.75$. On interpolation the
three recurrent branches lead the aggregate, yet the damped field DeepONet follows in fourth
place at $0.66$. It sits ahead of both latent one-shots and ahead of the \ac{mlp} baseline at
$0.94$, which holds the shedding. And the latent one-shot
control, which drops the autoregressive feedback, loses the oscillation entirely, $P = 0.03$
for the GRU branch and $0.00$ for the S4D branch.

Extrapolation lowers the held power to $P = 0.75$ for S4D and $0.60$ to $0.62$ for Mamba
and GRU and raises
$\rmsew$ for every model. The S4D branch reaches the encoder ceiling of
$\vKvsEncFloorExtrapP$ on this split, while the rollout \ac{rmse} of every model lies far above the
encoder floor of $\vKvsEncFloorExtrapRmse$. The Strouhal error stays low on this split even as amplitude and
structure degrade, so it is not on its own a quality verdict.

\subsection{Computational cost}
\label{sec:results-cost}

We apply the cost model of Sec.~\ref{sec:cost} to the regime champions. Solver wall-clock is measured
on a dual-socket AMD EPYC~9534 node ($16$ MPI ranks for the film, $4$ for the wake),
training and inference on a single NVIDIA RTX~PRO~6000 GPU, the two platforms stated separately.
Table~\ref{tab:cost} reports the per-case times, the training compute, the speedup $S$ of
Eq.~\eqref{eq:speedup}, and the break-even query count $Q^{*}$ of Eq.~\eqref{eq:breakeven}.

\begin{table*}[t]
\centering
\caption{Per-case cost of the regime champions. $t_{\mathrm{FEM}}$ on CPU, $t_{\mathrm{train}}$ and
$t_{\mathrm{infer}}$ on GPU, stated separately. $Q^{*}$ from Eq.~\eqref{eq:breakeven}. $^\ast$The \ac{cmp} solver
time is the measured as-run wall-clock at the $16$-core production point, which includes field output.
The compute-only time at $16$ cores is $9.2\,\si{\hour}$. The solver strong-scales further ($3.1\,\si{\hour}$
compute-only at $128$ cores, which would lower $S$ to about $1.2\times10^{4}$). We report the
configuration the dataset was produced with, which sits near the throughput optimum of the scaling
curve. The wake solver does not converge beyond $4$ ranks, so its production point is also its
fastest. $^\dagger$For the latent champions $t_{\mathrm{train}}$
is the convolutional-autoencoder fit ($\vCaeFitH\,\si{\hour}$) plus the dynamics fit ($\vSFourDFitH$ and $\vMambaFitH\,\si{\hour}$). The
encoder is shared by the latent wake models but charged in full to each champion.}
\label{tab:cost}
\footnotesize
\setlength{\tabcolsep}{3.5pt}
\begin{tabular}{@{}llccccc@{}}
\toprule
Regime & Champion & $t_{\mathrm{FEM}}$ & $t_{\mathrm{train}}$ & $t_{\mathrm{infer}}$ & $S$ & $Q^{*}$ \\
\midrule
\ac{cmp}$^\ast$ & field DeepONet & $11.3\,\si{\hour}$ & $0.43\,\si{\hour}$ & $0.91\,\si{\second}$ & $4.4\times10^{4}$ & $\approx 70$ \\
\ac{kvs}        & S4D-DON$^\dagger$    & $0.77\,\si{\hour}$ & $\vSFourDTrainH\,\si{\hour}$ & $0.66\,\si{\second}$ & $4.2\times10^{3}$ & $\approx \vKvsQstar$ \\
\ac{kvs}        & Mamba-DON$^\dagger$  & $0.77\,\si{\hour}$ & $\vMambaTrainH\,\si{\hour}$ & $0.60\,\si{\second}$ & $4.7\times10^{3}$ & $\approx \vKvsQstar$ \\
\bottomrule
\end{tabular}
\end{table*}

The per-query speedup is $\sim\!4.4\times10^{4}$ for the film and $\sim\!4.2$ to $4.7\times10^{3}$ for
the wake. The break-even query count is
$Q^{*} \approx 70$ for the film and $\approx \vKvsQstar$ for the wake, in each case close to the
$N_{\mathrm{train}}$ of the regime.
The training compute stays marginal. For each wake champion the shared convolutional-autoencoder
fit is charged in full, alongside its own dynamics fit, together $\vSFourDTrainH$ to
$\vMambaTrainH\,\si{\hour}$ or $\vKvsTrainQueriesLo$ to $\vKvsTrainQueriesHi$ solver queries. The film champion is a direct operator with no encoder stage, and its $0.43\,\si{\hour}$ amount
to $0.04$ queries. The surrogate therefore enters at its offline cost, $789\,\si{\hour}$ for the film and
$\vKvsOfflineH\,\si{\hour}$ for the wake, and undercuts the solver only beyond $Q^{*}$ queries. The implication for the
many-query setting is taken up in Sec.~\ref{sec:discussion}.

\section{Discussion}
\label{sec:discussion}

A surrogate validated on the simple regime does not stay reliable on the richer one, and the
winning architecture reverses between the regimes. The direct one-shot operator reconstructs the
film's process target at a relative \ac{cumwss} error of $0.027$, yet damps the wake oscillation it
should sustain to a shedding power of $P = 0.004$. The latent autoregressive models hold the wake,
with their structured state-space members at $P = 0.90$ and $0.91$, while on the film even the best of them,
the S4D branch, trails the direct operator by more than a factor of four on \ac{cumwss}, at $0.130$ against $0.027$.

The axis that decides this reversal is the treatment of time. Within each regime the latent models
share one frozen encoder, so between their one-shot and autoregressive variants the time
traversal changes while the representation and the branch family stay
fixed. The GRU branch occupies both of these cells in both regimes and carries the
simplest recurrence in the model set, so this pair isolates the time traversal without the
structured state-space machinery. On the wake the
autoregressive GRU holds the oscillation at $P = 0.90$ and the one-shot GRU damps it to $P = 0.03$,
and the S4D pair repeats the pattern at $0.90$ against $0.00$. On the film the order flips, and the
one-shot GRU reaches a \ac{cumwss} error of $0.099$ ahead of the autoregressive GRU at $0.167$.

Each regime rewards the time treatment that its dynamics demand. The wake is a self-sustained limit
cycle whose phase must be carried coherently across the full rollout. The autoregressive
feedback provides that memory, and a one-shot map does not, for a reason the latent
trajectories below make explicit. The film follows the prescribed
control ramp and develops no dynamics of its own, so a direct map from the boundary condition to the
field is sufficient and accumulates no rollout error. That a plain recurrent GRU already holds the
wake places the effect on the autoregressive feedback itself rather than on the structured
state-space parametrisation. Within the latent-AR family no branch leads across the board.
The phase-aligned \ac{rmse} puts the S4D and the Mamba marginally ahead of the GRU, the Strouhal
error puts the GRU ahead, and the shedding power sits at the encoder ceiling for all three. The
ranking within the family
therefore carries no architectural message. Under extrapolation the S4D branch holds the oscillation power of the frozen encoder's own
round-trip on these cases. The Mamba and GRU branches, at $0.60$ and $0.62$ against a ceiling
of $\vKvsEncFloorExtrapP$, lose a further part of it in their own dynamics. The rollout error
tells the opposite story. It lies far above the encoder floor for every branch, at $1.29$ to
$1.36$ against $\vKvsEncFloorExtrapRmse$, so the phase drift arises in the dynamics, not in the
representation.

\begin{figure*}[!t]
\centering
\includegraphics[width=\textwidth]{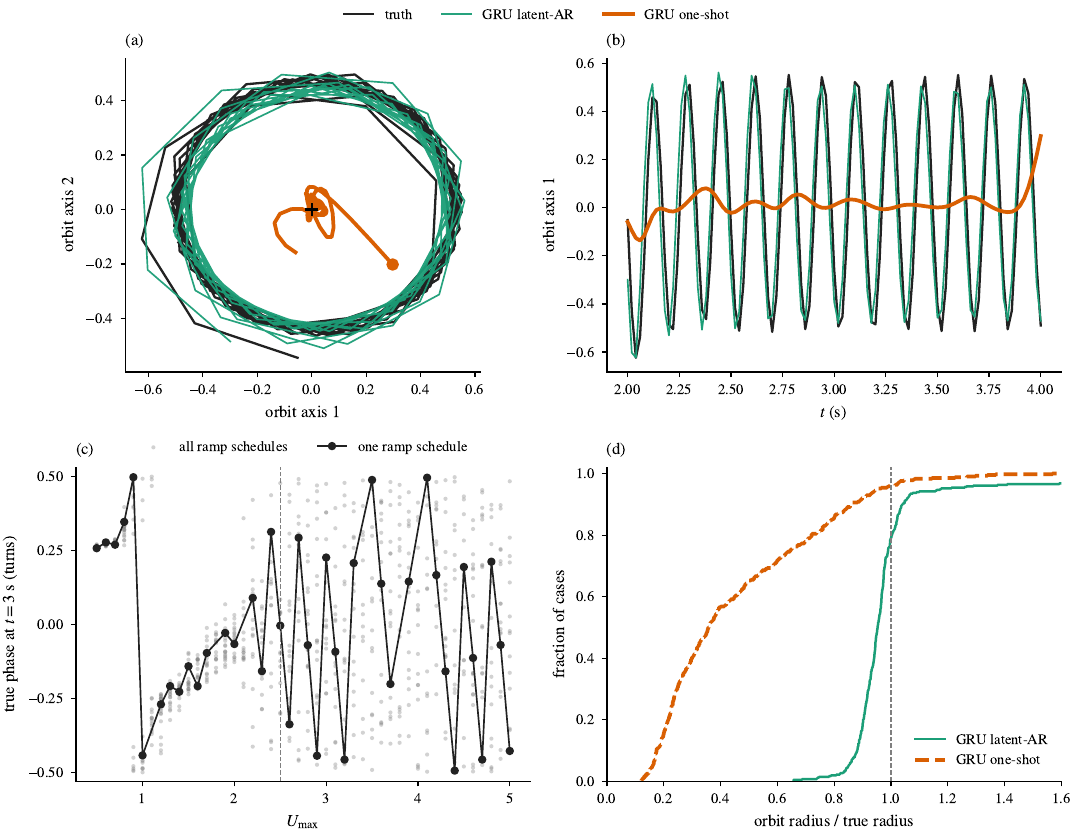}
\caption{Latent trajectories of the wake GRU pair on the interpolation split, latent-AR
against latent one-shot on the shared frozen encoder. (a)~Orbit of the test case with median
true amplitude, one seed, in the steady window $t \ge 2\,\si{\second}$, drawn in the plane of the two
leading directions of the mean-removed truth latents, the pair that carries the cycle, and
centred on the time mean of the truth. The cross
marks that mean, the dot the end of the one-shot trajectory. (b)~The first orbit direction over
time for the same case. (c)~Phase of the true cycle at $t = \vKvsOrbitTstar\,\si{\second}$, as a
fraction of one turn and wrapped to $\pm\tfrac12$, over $U_{\max}$ for all $\vKvsOrbitNall$ training and test cases (grey) and for one ramp schedule
(black, connected in the order of $U_{\max}$). The phase is the analytic-signal phase of the first
principal component of the truth latents. Below $U_{\max} \approx 1$ the true latent trace carries only a fraction of the
steady variance seen above it, so the wake has not developed shedding by this time and the phase
carries no meaning. The dashed line marks $U_{\max} = \vKvsOrbitUmaxSplit$.
(d)~Empirical cumulative distribution of the orbit radius relative to the truth over the $138$ test cases
and three seeds, in the late window $t \ge 2.5\,\si{\second}$ that excludes the initial transient. A
radius of one equals the true amplitude.}
\label{fig:kvs-latent-orbits}
\end{figure*}

The latent trajectories of the GRU pair show the mechanism directly
(Fig.~\ref{fig:kvs-latent-orbits}). In the steady window the
autoregressive GRU stays on the limit cycle of the truth. Its latent oscillation power is
$\vKvsOrbitPowerAr$ of the true one. Its orbit radius is $\vKvsOrbitRadiusAr$ in the median over cases and seeds, where one
equals the true amplitude. The one-shot GRU spirals into the centre of the cycle
(Fig.~\ref{fig:kvs-latent-orbits}a). It retains half of the latent variance, at
$\vKvsOrbitPowerOs$, but only $\vKvsOrbitBandOsPct\%$ of the variance in the shedding band of the
truth, against $\vKvsOrbitBandArPct\%$ for the autoregressive GRU. What remains is the initial
transient and a slow drift (Fig.~\ref{fig:kvs-latent-orbits}b), and once the transient has passed the orbit radius falls to
$\vKvsOrbitRadiusOsLate$ of the true one (Fig.~\ref{fig:kvs-latent-orbits}d). This is the
signature of a conditional mean. Fig.~\ref{fig:kvs-latent-orbits}c shows why the one-shot map
ends there. At a fixed time the phase of the true cycle is a function of the boundary parameters,
and the training ladder resolves that function only at low $U_{\max}$. Above
$U_{\max} = \vKvsOrbitUmaxSplit$ the shedding frequency has carried neighbouring cases more than a
quarter turn apart in $\vKvsOrbitAliasedHiPct\%$ of the neighbouring pairs, against
$\vKvsOrbitAliasedLoPct\%$ below. A smooth map from the parameters to the state cannot follow this
winding and averages the phase away, and the average over a cycle is its centre.
The winding sets a sampling bound. The $U_{\max}$ ladder of Sec.~\ref{sec:data} holds
$\vKvsLadderLevels$ levels at a step of $\vKvsLadderStep$. Since the shedding frequency rises
with $U_{\max}$, the phase difference between the top and the bottom of the ladder grows in
time. At $t = \vKvsOrbitTstar\,\si{\second}$ it amounts to $\vKvsOrbitWraps$ full turns, for
every ramp schedule alike. One step of the ladder thus shifts the phase by
$\vKvsLadderTurnsPerLevel$ of a turn. A map from the parameters can resolve the phase only if
one step shifts it by less than half a turn, the Nyquist condition applied to the parameter
axis~\cite{shannon1949communication}. That would need $\vKvsLadderNyquistLevels$ levels, or
$\vKvsLadderNyquistFactor$ times the present ladder, and the required number grows linearly with
the time since onset. Below this bound no map that is smooth in the parameters can recover the
phase, whatever its architecture. The autoregressive model needs no such map. It reads the phase
from its own state at every step and carries it forward, so the spacing of the parameter axis
never enters. Its error is the drift of that state, not an aliased average.

The S4D pair behaves the same way. The autoregressive S4D keeps $\vKvsOrbitSBandArPct\%$ of the
shedding-band variance at an orbit radius of $\vKvsOrbitSRadiusAr$. The one-shot S4D keeps
$\vKvsOrbitSBandOsPct\%$, at a latent power of $\vKvsOrbitSPowerOs$. In the plane of the cycle
both one-shot branches come to rest near the centre. After the transient their in-plane radius
is $\vKvsOrbitInPlaneOsLate$ and $\vKvsOrbitSInPlaneOsLate$, against
$\vKvsOrbitInPlaneTruthLate$ for the truth. The S4D one-shot fails in a
second way. Along the latent directions that the cycle does not use, its distance from the
truth's mean is $\vKvsOrbitSOutPlaneOsLate$ times the true orbit radius, against
$\vKvsOrbitOutPlaneTruthLate$ for the truth itself and $\vKvsOrbitOutPlaneOsLate$ for the GRU
one-shot. The S4D one-shot has thus averaged the phase away and, in addition, settled on a wrong
mean state.

The feedback loop thus has a benefit and a price, and the two regimes sit at its two ends. On
the wake the loop carries the phase, and without it the map returns the centre of the cycle. On
the film there is no phase to carry, and the loop only accumulates error along the rollout
(Fig.~\ref{fig:cmp-rollout-growth}). Between the two ends the benefit and the price would have
to be weighed. The two questions that close this section, whether the base state is stable and
how fast the field relaxes to its boundary condition, say which side dominates.

The representation matters too, but it does not decide the winner. On the film the direct field operator leads the best latent
one-shot variant by more than a factor of three on \ac{cumwss}, at $0.027$ against $0.099$, so mapping to the
field directly adds accuracy beyond the choice of time axis. On the wake no representation rescues
the one-shot map, and the direct and latent variants damp the shedding alike. Which architecture
wins a regime is set by the treatment of time. The representation changes the size of the
margin, a comparison the model set offers only within the one-shot pair, since the direct
autoregressive cell never trained stably, as Sec.~\ref{sec:limitations} reports.

The failure modes behind the reversal of the winning architecture are invisible to a single
aggregate error. On the wake's interpolation split the damped direct operator records a
rollout \ac{rmse} of $0.66$, inside the $0.33$ to $0.94$ range of the models that sustain the
shedding. The aggregate places it behind the three recurrent branches but ahead of the \ac{mlp}
baseline. The \ac{mlp} sustains the shedding at $P = 0.84$ and the direct operator holds none of
it at $P = 0.004$, yet the aggregate prefers the operator that has lost the oscillation over the
one that kept it. A field
relaxed toward its time mean stays pointwise close to a wake that
oscillates symmetrically about it. A surrogate that suppresses the shedding is therefore
barely penalised by the pointwise norm. On the extrapolation split the aggregate no longer hides
the failure. It inverts the ranking and puts the damped operator ahead of every model that sheds. On the
film the aggregate places the right model first but compresses its margin. The in-plane field error puts the direct
operator a factor of $\vFilmIpFactor$ ahead of the GRU-DON, the process target a factor of six.
In both regimes the
aggregate hides what matters, on the wake the lost oscillation, on the film the size of the
lead on the process target, and under wake extrapolation it inverts the ranking outright.
A selection by the aggregate alone would deploy the damped operator there.

Two further failure modes, pre-onset hallucination and onset latency, cut across the architecture
split rather than along it. On the interpolation split the GRU branch stays near rest before
onset at $H = 0.003$, on the same level as the damped field DeepONet. The S4D branch holds the
shedding just as well, yet moves more before onset than any model but the damped linear
baseline, at $H = 0.022$. Under extrapolation the state-space branches
lie at $H = 0.014$ to $0.017$, on the same level as the \ac{mlp} at $0.016$, and one seed in particular invents a full vortex street well before onset on individual cases (Fig.~\ref{fig:kvs-failures}b). The onset latency also cuts
across the split. The field DeepONet and the GRU one-shot reproduce the initial transient
on time yet never hold the oscillation, while among the shedders the \ac{mlp} is the latest by a
wide margin. Holding the oscillation therefore guarantees
neither a quiet pre-onset window nor a punctual onset, and the aggregate collapses these independent
axes into one number that rates each as a mild deviation.

\label{sec:z2}
On the wake's mirror symmetry the aggregate errs in the opposite direction. The wake admits two
mirror-image limit cycles, related by a reflection across the centreline and a half-period shift, and
under extrapolation nothing in the trained models pins which branch a prediction settles on. A prediction that develops
the correct dynamics on the mirror branch keeps the frequency, the shedding power, and the onset
unchanged, yet its pointwise rollout \ac{rmse} rises, since to a pointwise norm the mirror branch
reads as a phase shift of half a period, the largest possible on the limit cycle. The aggregate
that was too lenient with the
damped wake is therefore too harsh here, penalising an arbitrary branch choice as if it were a field
error. On the interpolation split the three recurrent models settle on the same branch as the
truth in all nine seeds, so the ambiguity is negligible in the aggregate there. Under extrapolation the
branch remains unpinned. In the mean over the three seeds the mirrored prediction is closer to the truth in $\vKvsMirrorCloserCases$
of the $\vKvsMirrorNCases$ cases. Scoring each case with whichever of the two is closer lowers the S4D
rollout \ac{rmse} by $\vKvsMirrorGainPct\%$. This closes most
of the gap to the damped operator. The
unpinned branch is also a deployment risk, which we record in
Sec.~\ref{sec:limitations}.

The practical rule for a \ac{cmp} surrogate programme follows from the mechanism. The choice
of surrogate should follow the dynamical character of the target flow rather than a validation
score from a simpler regime.
That character can be read off before any model is trained, from two questions about the
governing equations and the control schedule.

First, is the base state stable at the operating
point? The Stokes film is linear and dissipative, so every disturbance decays. The wake is
unstable above $\Reyn_c \approx 47$, which the design reaches at $U_{\max} \approx \vKvsUmaxCrit$.
It settles on a limit cycle whose phase is set by its own history, not by the momentary
boundary value. A direct map can carry that phase only if the parameter axis is sampled densely
enough to resolve it, which the ladder of this study does not. A model that feeds its own state
back needs no such sampling, and Fig.~\ref{fig:kvs-latent-orbits} shows what happens without
that feedback.

Second, how fast does the field relax to its boundary condition? The
film relaxes across the gap in $1.6\times10^{-3}\,\si{\second}$ against the one-second ramp, three orders of
magnitude apart, so the field is a function of the ramp position alone. The wake's shedding
period lies between $\vKvsShedPeriodMin$ and $\vKvsShedPeriodMax\,\si{\second}$ against the same one-second ramp window, so
the field depends on its history. A stable base state with a fast relaxation calls for a direct
map, an unstable one for autoregressive feedback. As a rule of thumb drawn from these two
regimes, we read the dynamical character as the ratio of the flow's own time scale, its
relaxation time or its oscillation period, to the time over which the boundary condition
changes. Far below one the flow follows its boundary
condition and a direct map suffices. Near or above one the flow carries its history and needs
feedback. The rule sorts our two regimes correctly.
Whether it holds beyond them is a question for further regimes, and it does not choose a model
within the autoregressive family.
Because neither a validation score from a simpler regime nor a single aggregate signals the
mismatch, the failure-mode-resolved metrics of Sec.~\ref{sec:metrics} are needed to expose it.
Extrapolation sharpens this caution. Under the extrapolation shift the ranking pattern of both regimes
persists. On the wake the held shedding power falls into the range $0.60$ to $0.75$. On the film the
direct operator's \ac{cumwss} error rises from $0.027$ to $0.032$. Several latent-AR film
metrics improve at the same time, so the direct operator's lead on the process target narrows
from a factor of $4.7$ to $3.7$.

\paragraph{The price of the speedup}
The speedup is the premise of a surrogate programme, not its result. Once trained, the regime
champions answer a query $4.4\times10^{4}$ times faster than the solver on the film and $4.2$ to
$4.7\times10^{3}$ times faster on the wake (Sec.~\ref{sec:results-cost}). The gain materialises only
beyond a break-even query count that the training set dominates, $Q^{*} \approx 70$ for the film and
$\approx \vKvsQstar$ for the wake, because the offline cost is carried by the $N_{\mathrm{train}}$ reference
simulations rather than by the training compute. The speedup therefore pays off only in
many-query settings, such as parameter studies, inverse problems, and uncertainty
quantification, and a single prediction gains nothing. Reporting the raw speedup
without the break-even count would overstate the practical gain, a reporting pattern that a recent
survey identifies as a source of overoptimism in machine learning for
\acp{pde}~\cite{mcgreivy2024weak}.
\section{Limitations and Outlook}
\label{sec:limitations}

Several limitations bound these conclusions, and each closes with the outlook it motivates.

The first limitation is the comparison design. The two regimes differ in encoder,
dimensionality, dynamics, and the instantiation of the ramp family (Sec.~\ref{sec:data})
at once, so encoder and regime are confounded. Apart from the
within-wake time-treatment ablation (Sec.~\ref{sec:discussion}), the reading that the regime
dictates the architecture is a consistent observation across two regimes, not a controlled
result~\cite{lu2022faircomparison}. A named best model is the strongest entry in this comparison,
not a certified optimum.

The second is the physical scope. The wake data are a two-dimensional direct numerical
simulation, an idealisation that ignores the three-dimensional effects a real wake develops at higher
Reynolds numbers. This affects the upper end of our range, including the whole
extrapolation split. The film covers the hydrodynamic
sub-problem on a fixed geometry, with the headline \ac{cumwss} a removal-relevant but partial
target (Sec.~\ref{sec:data}). A more detailed process simulation is the
application-facing next step, for example a particle-transport model of material removal.

The third concerns the reference itself. The surrogate error is measured against a
finite-element solution whose discretisation error we have not assessed, so the reported numbers
quantify agreement with the reference rather than with the physics. A mesh and time-step
convergence study could bound the discretisation error of the reference.

The generalisation we establish is bounded to the trained input family. The controls vary only
within the multi-stage piecewise-linear ramp family, so what we validate is an in-dis\-tribution
emulator over this family, not an operator over an arbitrary space of controls. Within the
family the transient remains the hardest part, and conditioning the rollout on the onset is the
natural follow-up.

The encoders are frozen and optimised for reconstruction rather than for predictability, so the
representation itself can cap what the latent pipelines can reach. On the wake it caps the
amplitude. The encoder alone returns a shedding power of $\vKvsEncFloorInterpP$ on interpolation
and $\vKvsEncFloorExtrapP$ on extrapolation, and the best recurrent branch reaches that ceiling
on both splits within the seed spread, while the rollout error of every branch stays far above the encoder floor. The
cap works against the latent models only. They win the wake at the ceiling, so the wake ranking
does not depend on it. On the film the truncated \ac{pod} basis is part of the latent design. Its
reconstruction floor on the process target, $\vFilmEncFloorInterpCw$ against $0.027$ for the
direct operator on interpolation, lies above the direct operator's error, so the representation
carries part of the film margin, as Sec.~\ref{sec:discussion} states. It does not carry the loss.
On the bulk field the floor lies below the direct operator, at $\vFilmEncFloorInterpIp$ against
$0.23$, and the latent rows still lose there. Richer bases lower the floor further but did not
improve the trained rows in our trials. With more state to carry, the rollouts either diverged or
converged to worse values. Encoders trained jointly
with the dynamics~\cite{lusch2018koopman}, or parameter-adaptive
bases~\cite{franco2024dod}, could lower it.

The configurations were found by exploration under
the tuning asymmetry stated in Sec.~\ref{sec:training}, not by a cross-validated search over
every cell. A pre-registered, budget-matched search is the rigorous follow-up. A post-hoc
search over the film's latent autoregressive rows, including the rollout horizon, saturated at
the \ac{pod} floor and left the ranking unchanged. A full search would sharpen the margins but
not change their sign. The flip rests on the two mechanisms
of Sec.~\ref{sec:discussion}. A one-shot map collapses to the conditional mean where the phase is
aliased across cases, and an autoregressive rollout accumulates error where nothing is to be
carried. Tuning removes neither. The first is a limit of the data, not of the architecture. At
$t = \vKvsOrbitTstar\,\si{\second}$ the $U_{\max}$ ladder would need $\vKvsLadderNyquistFactor$
times as many levels to reach the sampling bound of Sec.~\ref{sec:discussion}. The later the
time, the more levels it needs. Whether a one-shot map then learns the phase is untested.

The surrogate also
yields point predictions without error certificates, which certified reduced-order modelling
provides~\cite{benner2015survey} and operator learning still seeks~\cite{qiu2025variational}. The
screen-then-verify workflow is therefore motivated rather than guaranteed. A decision rule could
come from a physics-residual indicator on the predicted field, provided the residual tracks the
true error, or from an eigenvalue analysis of the learned autoregressive step
map~\cite{mohan2024illusion}.

The model set leaves the direct autoregressive cell of Table~\ref{tab:models}
unfilled. It is the one cell in which the claim that no single architecture wins both
regimes could be overturned. The mechanism predicts what a model there would do: carry the
feedback that holds the wake, and accumulate the rollout error that costs every autoregressive
model the film. We trained the field DeepONet autoregressively, with a grown rollout horizon and
noise injection, but the rollout never stabilised. We also tested an autoregressive Fourier
neural operator \cite{li2021fno} on the rasterised wake, but did not get it to work at practical training
cost. On the thin three-dimensional \ac{cmp} film it would additionally require a
geometry-adapted variant~\cite{li2023gino,peyvan2026fusion}. The cell remains open for future
attempts.

The statistical claims are bounded by the sampling. Each stochastic model is the mean over
three seeds, and we do not read a difference as a ranking where it falls within the seed spread.
The wake extrapolation split is underpowered, its $200$ cases sharing only ten distinct
$U_{\max}$ levels, so several differences among its leading models fall within the seed spread.
The film splits hold $15$ test cases each. Cross-validation over repeated splits would add the
split variance to these bounds, at the cost of retraining the full model matrix per fold.

The training-set size is itself a design variable.
Because the break-even query count of Sec.~\ref{sec:results-cost} is set by
$N_{\mathrm{train}}$, the accuracy attainable at smaller training sets determines how far it can
drop. A first sweep over reduced training sets, on a single subset draw, points to saturation
well below the full training-set sizes in both regimes, a direction rather than a result.
Characterising the curve with disjoint draws and confidence intervals is the next measurement on
this axis, where acquisition strategies such as active learning could lower the training-set
size that sets the break-even~\cite{musekamp2024al4pde}.

A final note concerns the wake's unpinned mirror branch (Sec.~\ref{sec:z2}). This is an
observation rather than a limitation of the study, but under extrapolation it leaves a deployed model free to settle
on the mirrored limit cycle. On the interpolation split every seed settles on the same
branch as the truth (Sec.~\ref{sec:z2}). Symmetry-aware modelling, for example an equivariant
parametrisation of the latent dynamics, could pin the branch and has improved rollout stability in
related settings~\cite{huang2025geometric}.

\FloatBarrier
\section{Conclusion}
\label{sec:conclusion}

Motivated by the cost of physical \ac{cmp} experimentation, we asked whether the surrogate
architecture validated on a simple regime stays reliable when it is trained from scratch on richer
physics. A two-regime
contrast, a linear Stokes film against an oscillatory \ac{kvs} wake, answers it on one shared
evaluation pipeline. The answer is a no-free-lunch result for the architectures tested, at a
given budget of training simulations. A
one-shot direct operator wins the boundary-driven film, latent autoregressive models win the
self-sustained wake, and the deciding axis is the treatment of time. The wake requires the phase
memory that autoregressive feedback provides, while the film rewards the direct map that
accumulates no rollout error.

Two measurement findings frame this result. A single aggregate error hides the damped
oscillation on the wake, compresses the sixfold lead on the film, and inverts the model ranking
under wake extrapolation, so failure modes must be
read from the suite of five physical questions, accuracy, structure, invented motion, amplitude,
and timing. And the per-query speedup of three to four orders of magnitude pays off only beyond
a break-even query count set almost entirely by the size of the training set, so the surrogate
is a many-query instrument.

The practical implication for \ac{cmp} surrogate programmes is that cheap validation on a
simplified regime is not evidence of an architecture's reliability on a richer one. The choice of surrogate should
follow the dynamical character of the target flow, read as the ratio of the flow's own time
scale to that of its boundary condition, and its verification a failure-mode-resolved
suite. Neither the contrast design nor the evaluation methodology is tied to the architectures
compared here.

\section*{CRediT authorship contribution statement}
\textbf{Georg Winkler:} Conceptualization, Methodology, Software, Formal analysis, Investigation,
Data curation, Validation, Visualization, Writing -- original draft, Writing -- review \& editing.
\textbf{Martin Stoll:} Conceptualization, Resources, Supervision, Writing -- review \& editing, Funding acquisition.

\section*{Declaration of competing interest}
The authors declare that they have no known competing financial interests or personal relationships
that could have appeared to influence the work reported in this paper.

\section*{Data availability}
The configuration files, evaluation harness, value layer behind every printed number, and metric definitions are versioned in a
repository that will be released with the publication. The raw simulation datasets and trained checkpoints are archived on institutional scratch
storage and available from the authors on reasonable request.

\section*{Acknowledgment}
We sincerely thank the European Social Fund (ESF) and the Free State of Saxony
of the Federal Republic of Germany for funding this work under grant number 100693458 (project
WiProFlex).
\par\medskip\noindent
{\centering\includegraphics[width=0.8\columnwidth]{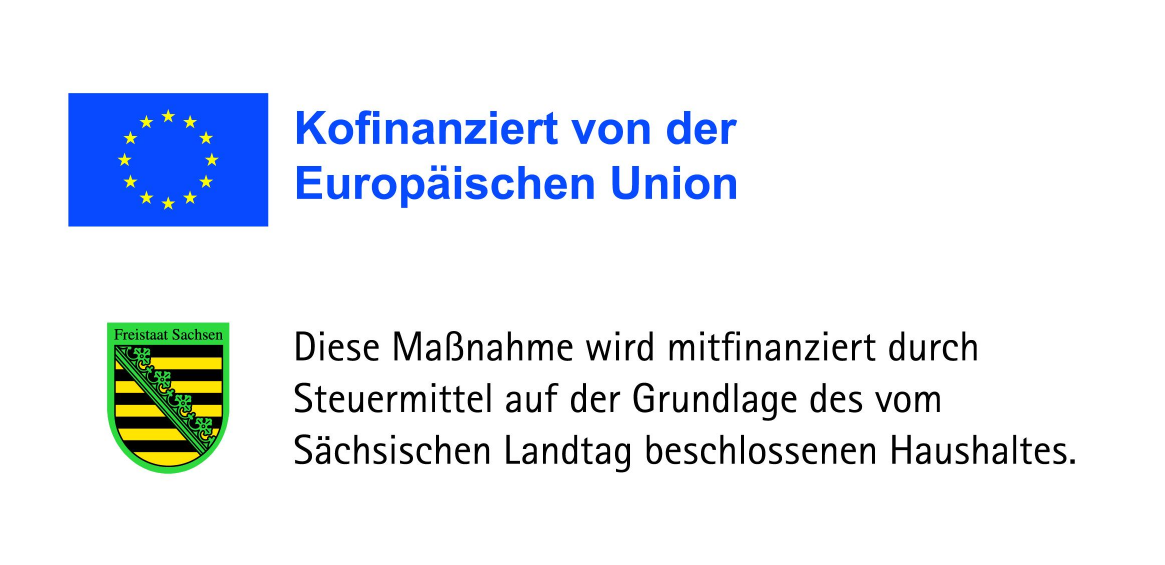}\par}

\section*{Declaration of generative AI and AI-assisted technologies in the writing process}
While preparing this work, the authors used Anthropic's large language model Claude to assist in
drafting and revising the text and in writing and reviewing the research software and
data-analysis scripts. The authors designed the study, carefully reviewed and verified all
content and results, and take full responsibility for the final version of the manuscript.

\end{document}